\documentclass{article}
\usepackage{amsfonts}
\usepackage{amsmath}
\usepackage{amssymb}
\usepackage{geometry}
\usepackage{setspace}
\usepackage{sectsty}
\usepackage{lineno}
\usepackage{theorem}
\usepackage{fancyhdr}
\usepackage{graphicx}

\newtheorem{theorem}{Theorem}

\newtheorem{definition}[theorem]{Definition}

\newtheorem{lemma}[theorem]{Lemma}

\newtheorem{remark}[theorem]{Remark}

\newenvironment{proof}[1][Proof]{\noindent\textbf{#1.} }{\ \rule{0.5em}{0.5em}}
\DeclareMathOperator{\pen}{pen}
\DeclareMathOperator{\crit}{crit}
\DeclareMathOperator{\leb}{Leb}
\DeclareMathOperator{\card}{Card}

\DeclareMathOperator{\Span}{Span}

\DeclareMathOperator{\p}{p}
\input{tcilatex}

\begin{document}

\title{Optimal model selection in heteroscedastic regression using piecewise
polynomial functions}
\author{A. Saumard\thanks{%
Research partly supported by the french Agence Nationale de la Recherche
(ANR 2011 BS01 010 01 projet Calibration), NI-AID grant 2R01 AI29168-04 and
a PIMS postdoctoral fellowship.} \\
Department of Statistics, University of Washington, Seattle, WA 98195, USA\\
INRIA Saclay \^{I}le-de-France, France}
\date{April 4, 2013}
\maketitle

\begin{abstract}
We consider the estimation of a regression function with random design and
heteroscedastic noise in a nonparametric setting. More precisely, we address
the problem of characterizing the optimal penalty when the regression
function is estimated by using a penalized least-squares model selection
method. In this context, we show the existence of a minimal penalty, defined
to be the maximum level of penalization under which the model selection
procedure totally misbehaves. The optimal penalty is shown to be twice the
minimal one and to satisfy a non-asymptotic pathwise oracle inequality with
leading constant almost one. Finally, the ideal penalty being unknown in
general, we propose a hold-out penalization procedure and show that the
latter is asymptotically optimal.\bigskip

\end{abstract}

\section{Introduction}

Given a collection of models and associated estimators, two different model
selection tasks can be tackled: find out the smallest true model
(consistency problem), or select an estimator achieving the best performance
according to some criterion, called a risk or a loss\textit{\ }(efficiency
problem). We focus on the efficiency problem, where the leading idea of
penalization, that goes back to early works of Akaike 
\cite{Akaike:70,Akaike:73}
and Mallows \cite{Mallows:73}, is to perform an unbiased - or uniformly
biased - estimation of the risk of the estimators. FPE and AIC procedures
proposed by Akaike respectively in \cite{Akaike:70} and \cite{Akaike:73}, as
well as Mallows' $C_{p}$ or $C_{L}$ \cite{Mallows:73}, aim to do so by
adding to the empirical risk a penalty which depends on the dimension of the
models.

The first analysis of such procedures had the drawback of being
fundamentally asymptotic, considering in particular that the number of
models as well as their dimensions are fixed while the sample size tends to
infinity. As explained for instance in Massart \cite{Massart:07}, in various
statistical settings it is natural to let these quantities depend on the
amount of data. Thus, pointing out the importance of Talagrand's type
concentration inequalities in the nonasymptotic approach, Birg\'{e} and
Massart 
\cite{BirMassart:97,BirgeMassart:01}
and Barron, Birg\'{e} and Massart \cite{BarBirMassart:99} have been able to
build nonasymptotic oracle inequalities for penalization procedures. Their
framework takes into account the complexity of the collection of models as a
parameter depending on the sample size.

In an abstract risk minimization framework, which includes statistical
learning problems such as classification or regression, many
distribution-dependent and data-dependent penalties have been proposed, from
the more general and less accurate global penalties, see Koltchinskii \cite%
{Kolt:01}, Bartlett et al. \cite{BarBouLugosi:02}, to the refined local
Rademacher complexities in the case where some favorable noise conditions
hold (see for instance Bartlett, Bousquet and Mendelson \cite%
{BarBouMendelson:05}, Koltchinskii \cite{Kolt:06}). But as a price to pay
for generality, the above penalties suffer from their dependence on unknown
constants. These penalized procedures are very difficult to implement and
calibrate in practice. Moreover, the existing risk bounds for these
procedures contain very large leading constants. Other general-purpose
penalties have been proposed, such as the bootstrap penalties of Efron \cite%
{Efron:83} and the resampling and $V$-fold penalties of Arlot 
\cite{Arl:2008a,Arl:09}%
. These penalties are essentially resampling estimates of the difference
between the empirical risk and the risk. Arlot 
\cite{Arl:2008a,Arl:09}
proved sharp pathwise oracle inequalities for the resampling and $V$-fold
penalties in the case of regression with random design and heteroscedastic
noise on histograms models, and conjectured that the restriction to
histograms is mainly technical and that his results can be extended to more
general situations.

Model selection \textit{via} penalization is not the only method which
provides sharp oracle inequalities for the estimation of a nonparametric
regression function. Indeed, aggregation techniques and PAC-Bayesian bounds
also allow to obtain nearly optimal constants in the oracle inequalities.
Bunea et al.\textit{\ }\cite{BunTsyWeg:07} derived some sharp oracle
inequalities for different aggregation tasks by means of a single unifying
procedure. However, the authors asked for a fixed design and homoscedastic
Gaussian noise. By using aggregation with exponential weights, Dalalyan and
Tsybakov obtained in \cite{DalTsy:07} oracle inequalities of a PAC-Bayesian
flavor with leading constant one and optimal rate of the remainder term for
the estimation of a regression function with deterministic design and
homoscedastic errors. Furthermore, these authors allowed error distributions
which are symmetric or $n$-divisible. PAC-Bayesian methods are
systematically investigated in Catoni, \cite{Catoni:04}. The work of Lecu%
\'{e} and Mendelson \cite{LecMen:09} concerning the aggregation by empirical
risk minimization of a finite family of functions seems to handle the case
of a random design and heteroscedastic noise, even if this example is not
explicitly developed. The oracle inequalities obtained by Lecu\'{e} and
Mendelson are sharp and valid with probability close to one. In particular,
they are related to oracle inequalities obtained, in expectation, by Catoni
in \cite{Catoni:04}.

A difference between aggregation and model selection studies, is that in
most aggregation results, the estimators at hand are considered as
deterministic functions. However, notable exceptions are the following.
Leung and Barron \cite{LeungBarron:06} proved sharp oracle inequalities for
the aggregation of projection estimators in the Gaussian sequence model.
Rigollet and Tsybakov \cite{RigTsy:12} recently showed sharp bounds for the
aggregation of some linear estimators, including projection estimators, in a
regression setting, with fixed design and homoscedastic Gaussian noise. More
general PAC-Bayesian type inequalities were also recently obtained by
Dalalyan and Salmon \cite{DalSal:12}, considering the aggregation of affine
estimators in heteroscedastic regression, with Gaussian noise and fixed
design.

Birg\'{e} and Massart \cite{BirMas:07} discovered, in a generalized linear
Gaussian model setting, that the optimal penalty is closely related to the
minimal one. An optimal penalty is a penalty which gives an oracle
inequality with leading constant converging to one when the sample size
tends to infinity. The minimal penalty is defined to be the maximal penalty
under which the procedure totally misbehaves (in a sense to be specified
below). Birg\'{e} and Massart \cite{BirMas:07} proved sharp upper and lower
bounds for the minimal penalty. These authors also showed that the optimal
penalty is twice the minimal one, both for small and large collections of
models. These facts are called the\textit{\ slope heuristics}. The authors
also exhibited a jump in the dimension of the selected model occurring
around the value of the minimal penalty, and used it to estimate the minimal
penalty from the data. Taking a penalty equal to twice the previous estimate
then gives a nonasymptotic quasi-optimal data-driven model selection
procedure. The algorithm proposed by Birg\'{e} and Massart \cite{BirMas:07}
to estimate the minimal penalty relies on the previous knowledge of the
shape of the latter, which is a known function of the dimension of the
models in their setting. Thus, their procedure gives a data-driven \textit{%
calibration }of the minimal penalty.

Considering the case of Gaussian least-squares regression with unknown
variance, Baraud et al. \cite{BarGirHuet:09} have also derived lower bounds
on the penalty terms for small and large collections of models. In the
setting of maximum likelihood estimation of density on histograms, Castellan 
\cite{Castellan:99} obtained a lower bound on the penalty term, in the case
of small collections of models.

The slope heuristics has been then extended by Arlot and Massart \cite%
{ArlotMassart:09} in a bounded regression framework, with heteroscedastic
noise and random design. The authors considered least-squares estimators on
a \textquotedblleft small\textquotedblright\ collection of histograms
models. Their analysis differs from the one of Birg\'{e} and Massart \cite%
{BirMas:07} in an important way. Indeed, Arlot and Massart \cite%
{ArlotMassart:09} did not assume a particular \textit{shape} of the penalty
term. As a matter of fact, the penalties considered by Birg\'{e} and Massart 
\cite{BirMas:07} were known functions of the dimension of the models,
whereas heteroscedasticity of the noise allowed Arlot and Massart to
consider situations where the shape of the penalty is not even a function of
the dimension of the models. In such general cases, the authors proposed to
estimate the shape of the penalty by using Arlot's resampling or $V$-fold
penalties, proved to be efficient in their regression framework by Arlot 
\cite{Arl:2008a,Arl:09}%
.

The approach developed in \cite{ArlotMassart:09} is more general than the
histogram case, except for some identified technical parts of the proofs,
thus providing a general framework that can be applied to other problems.
The authors have also identified, in the case of histograms, the minimal
penalty as the mean of the empirical excess loss on each model, and the
ideal penalty to be estimated as the sum of the empirical excess loss and
true excess loss on each model. The slope heuristics then heavily relies on
the fact that the empirical excess loss is equivalent to the true excess
loss for models of reasonable dimensions.

Arlot and Massart \cite{ArlotMassart:09} conjectured that this equivalence
between the empirical and true excess loss is a quite general fact in
M-estimation. A general result supporting this conjecture is the high
dimensional Wilks' phenomenon investigated by Boucheron and Massart \cite%
{BouMas:10} in the setting of bounded contrast minimization. The authors
derive in \cite{BouMas:10} concentration inequalities for the empirical
excess loss, under some margin conditions (called \textquotedblleft noise
conditions\textquotedblright\ by the authors) and when the considered model
satisfies some general \textquotedblleft complexity
condition\textquotedblright\ on the first moment of the supremum of the
empirical process on localized slices of variance in the loss class. The
latter assumption can be explicated under suitable covering entropy
conditions on the model.

Lerasle \cite{Ler:12} proved the validity of the slope heuristics in a
least-squares density estimation setting, under rather mild conditions on
the considered linear models. The approach developed by the author in this
framework allows sharp computations and the empirical excess loss is shown
to be exactly equal to the true excess loss. Lerasle \cite{Ler:12} also
proved in the least-squares density estimation setting the efficiency of
Arlot's resampling penalties. Moreover, Lerasle \cite{Ler:11} generalized
the previous results to weakly dependent data. Arlot and Bach \cite%
{Arl_Bac:2009} recently considered the problem of selecting among linear
estimators in nonparametric regression. Their framework includes model
selection for linear regression, the choice of a regularization parameter in
kernel ridge regression or spline smoothing, and the choice of a kernel in
multiple kernel learning. In such cases, the minimal penalty is not
necessarily half the optimal one, but the authors propose to estimate the
unknown variance by the minimal penalty and to use it in a plug-in version
of Mallows' $C_{L}$. The latter penalty is proved to be optimal by
establishing a nonasymptotic oracle inequality with constant close to one,
converging to one when the sample size tends to infinity.

In this paper, we prove the validity of the slope heuristics in the
framework of bounded regression with random design and heteroscedastic
noise. This is done by considering a \textquotedblleft
small\textquotedblright\ collection of finite-dimensional linear models of
piecewise polynomial functions. This setting extends the case of histograms
already treated by Arlot and Massart \cite{ArlotMassart:09}. An interesting
consequence is that piecewise polynomial functions are known to have good
approximation properties in Besov spaces and can lead to minimax rates of
convergence, see for instance 
\cite{BarBirMassart:99,Tsy:04b}%
. As a matter of fact, histograms allow minimax procedures only on H\"{o}%
lder spaces.

Our validation of the slope heuristics is of asymptotic nature. However, the
complexity of the collection of models as well as their dimensions are not
constant terms in our analysis. These quantities are indeed allowed to
depend on the sample size $n$.

If the noise is homoscedastic, then the shape of the ideal penalty is known,
and is linear in the dimension of the models as in the case of Mallows' $%
C_{p}$. However, if the noise is heteroscedastic, then Arlot \cite{Arl:2010}
showed that the ideal penalty is not even a function of the linear
dimensions of the models. So, it is necessary to give a suitable estimator
of this shape. As emphasized by Arlot 
\cite{Arl:2008a,Arl:09}%
, $V$-fold and resampling penalties are good, natural candidates for this
task. In this paper, we show that a hold-out penalty - which is closely
related to a special case of resampling penalty - is indeed asymptotically
optimal under very mild conditions on the data split. As a matter of fact, a
half-and-half split leads to an optimal penalization. It is worth noticing
that hold-out type procedures have also been exploited in Chapter 8 of
Massart \cite{Massart:07} as simple tools to overcome the margin adaptivity
issue in classification.

The paper is organized as follows. In Section \ref{section_framework_reg
copy(1)}, we describe the statistical framework. The slope heuristics is
presented in Section \ref{section_slope_heuristics}, and the hold-out
penalization is considered in Section \ref{section_hold_out}. The proofs are
collected in Section \ref{section_proof_slope_reg}.

\section{Statistical framework\label{section_framework_reg copy(1)}}

\subsection{Penalized least-squares model selection}

Let us take $n$ independent observations $\xi _{i}=\left( X_{i},Y_{i}\right)
\in \mathcal{X\times }\mathbb{R}$ with common distribution $P$. In Sections %
\ref{section_piecewise_polynomials}, \ref{section_main_assumptions}-\ref%
{section_hold_out} the feature space $\mathcal{X=}\left[ 0,1\right] $. The
marginal distribution of $X_{i}$ is denoted by $P^{X}$. We assume that the
data satisfy the following relation%
\begin{equation}
Y_{i}=s_{\ast }\left( X_{i}\right) +\sigma \left( X_{i}\right) \varepsilon
_{i}\text{ },  \label{regression_model_2}
\end{equation}%
where $s_{\ast }\in L_{2}\left( P^{X}\right) $. Conditionally to $X_{i}$,
the residual $\varepsilon _{i}$ is assumed to have zero mean and variance
equal to one. The function $\sigma :\mathcal{X\rightarrow }\mathbb{R}_{+}$
is the unknown heteroscedastic noise level. A generic random variable with
distribution $P$, independent of the sample $\left( \xi _{1},...,\xi
_{n}\right) $, is denoted by $\xi =\left( X,Y\right) $.

It follows from (\ref{regression_model_2}) that $s_{\ast }$ is the unknown
regression function of $Y$ with respect to $X$. Our aim is to estimate $%
s_{\ast }$ from the sample. To do so, we are given a finite collection of
models $\mathcal{M}_{n}$, with cardinality depending on the sample size $n$.
Each model $M\in \mathcal{M}_{n}$ is assumed to be a finite-dimensional
vector space. We denote by $D_{M}$ the linear dimension of $M$. In the main
part of this paper, we focus on models of piecewise polynomial functions,
that are introduced in Section \ref{section_piecewise_polynomials} below.

We denote by $\left\Vert s\right\Vert _{2}=\left( \int_{\mathcal{X}%
}s^{2}dP^{X}\right) ^{1/2}$ the usual norm in $L_{2}\left( P^{X}\right) $
and by $s_{M}$ the linear projection of $s_{\ast }$ onto $M$ in the Hilbert
space $\left( L^{2}\left( P^{X}\right) ,\left\Vert \cdot \right\Vert
_{2}\right) $. For a function $f\in L_{1}\left( P\right) $, we write $%
P(f)=Pf=\mathbb{E}\left[ f\left( \xi \right) \right] $. By setting $%
K:L_{2}\left( P^{X}\right) \rightarrow L_{1}\left( P\right) $ the
least-squares contrast, defined by 
\begin{equation}
K\left( s\right) :\left( x,y\right) \mapsto \left( y-s\left( x\right)
\right) ^{2}\text{ , \ \ \ \ }s\in L_{2}\left( P^{X}\right) \text{ ,}
\label{def_contrast}
\end{equation}%
the regression function $s_{\ast }$ satisfies 
\begin{equation}
s_{\ast }=\arg \min_{s\in L_{2}\left( P^{X}\right) }P\left( K\left( s\right)
\right) \text{ .}  \label{def_target}
\end{equation}%
For the linear projections $s_{M}$ we get 
\begin{equation}
s_{M}=\arg \min_{s\in M}P\left( K\left( s\right) \right) \text{ .}
\label{def_projection}
\end{equation}%
For each model $M\in \mathcal{M}_{n}$, we consider a least-squares estimator 
$s_{n}\left( M\right) $ (possibly non unique), satisfying 
\begin{align*}
s_{n}\left( M\right) & \in \arg \min_{s\in M}\left\{ P_{n}\left( K\left(
s\right) \right) \right\} \\
& =\arg \min_{s\in M}\left\{ \frac{1}{n}\sum_{i=1}^{n}\left( Y_{i}-s\left(
X_{i}\right) \right) ^{2}\right\} \text{ ,}
\end{align*}%
where $P_{n}=n^{-1}\sum_{i=1}^{n}\delta _{\xi _{i}}$ is the empirical
measure built from the data.

In order to avoid cumbersome notations, we will often write $Ks$ in place of 
$K\left( s\right) $ for the image of a suitable function $s$ by the contrast 
$K$. We measure the performance of the least-squares estimators by their
excess loss,%
\begin{equation*}
\ell \left( s_{\ast },s_{n}\left( M\right) \right) :=P\left( Ks_{n}\left(
M\right) -Ks_{\ast }\right) =\left\Vert s_{n}\left( M\right) -s_{\ast
}\right\Vert _{2}^{2}\text{ .}
\end{equation*}%
We have the following decomposition,%
\begin{equation*}
\ell \left( s_{\ast },s_{n}\left( M\right) \right) =\ell \left( s_{\ast
},s_{M}\right) +\ell \left( s_{M},s_{n}\left( M\right) \right) \text{ ,}
\end{equation*}%
where%
\begin{equation*}
\ell \left( s_{\ast },s_{M}\right) :=P\left( Ks_{M}-Ks_{\ast }\right)
=\left\Vert s_{M}-s_{\ast }\right\Vert _{2}^{2}\text{ \ \ \ and \ \ \ }\ell
\left( s_{M},s_{n}\left( M\right) \right) :=P\left( Ks_{n}\left( M\right)
-Ks_{M}\right) \geq 0\text{ .}
\end{equation*}%
The quantity $\ell \left( s_{\ast },s_{M}\right) $ is called the bias of the
model $M$ and $\ell \left( s_{M},s_{n}\left( M\right) \right) $ is the
excess loss of the least-squares estimator $s_{n}\left( M\right) $ on the
model $M$. By the Pythagorean identity, we have%
\begin{equation*}
\ell \left( s_{M},s_{n}\left( M\right) \right) =\left\Vert s_{n}\left(
M\right) -s_{M}\right\Vert _{2}^{2}\text{ .}
\end{equation*}

Given the collection of models $\mathcal{M}_{n}$, an oracle model $M_{\ast }$
is defined as a minimizer of the losses - or equivalently excess losses - of
the estimators at hand,%
\begin{equation}
M_{\ast }\in \arg \min_{M\in \mathcal{M}_{n}}\left\{ \ell \left( s_{\ast
},s_{n}\left( M\right) \right) \right\} \text{ .}  \label{oracle_model}
\end{equation}%
The associated oracle estimator $s_{n}\left( M_{\ast }\right) $ thus
achieves the best performance in terms of excess loss among the collection $%
\left\{ s_{n}\left( M\right) ;M\in \mathcal{M}_{n}\right\} $. The oracle
model is a random quantity because it depends on the data and it is also
unknown as it depends on the distribution $P$ of the data. We propose to
estimate the oracle model by a penalization procedure.

Given some known penalty $%
\pen%
$, that is a function from $\mathcal{M}_{n}$ to $\mathbb{R}$, we consider
the following data-dependent model, also called selected model,%
\begin{equation}
\widehat{M}%
\in \arg \min_{M\in \mathcal{M}_{n}}\left\{ P_{n}\left( Ks_{n}\left(
M\right) \right) +%
\pen%
\left( M\right) \right\} \text{ }.  \label{def_proc_2_reg}
\end{equation}%
Our aim is then to find a good penalty, such that the selected model $%
\widehat{M}%
$ satisfies an oracle inequality of the form%
\begin{equation*}
\ell \left( s_{\ast },s_{n}\left( 
\widehat{M}%
\right) \right) \leq C\times \ell \left( s_{\ast },s_{n}\left( M_{\ast
}\right) \right) \text{ ,}
\end{equation*}%
with some positive constant $C$ as close to one as possible and with
probability close to one, typically more than $1-Ln^{-2}$ for some positive
constant $L$.

\subsection{Piecewise polynomial functions\label%
{section_piecewise_polynomials}}

Let us take $\mathcal{X=}\left[ 0,1\right] $ the unit interval and $\mathcal{%
P}$ a finite partition of $\mathcal{X}$. For a positive integer $r$ and any $%
\left( I,j\right) \in \mathcal{P\times }\left\{ 0,...,r\right\} $, we set 
\begin{equation*}
p_{I,j}:x\in \mathcal{X}\mapsto x^{j}\mathbf{1}_{I}\left( x\right) \text{ .}
\end{equation*}

\begin{definition}
\label{def_model_piece_pol}A finite dimensional vector space $M$ is said to
be a model of piecewise polynomial functions, with respect to the finite
partition $\mathcal{P}$ of $\mathcal{X=}\left[ 0,1\right] $ and of degrees
not larger than $r\in \mathbb{N}$, if 
\begin{equation*}
M=%
\Span%
\left\{ p_{I,j}\text{ };\text{ }\left( I,j\right) \in \mathcal{P\times }%
\left\{ 0,...,r\right\} \right\} \text{ .}
\end{equation*}%
The linear dimension of $M$ is then equal to $\left( r+1\right) \left\vert 
\mathcal{P}\right\vert $.
\end{definition}

Notice that models of histograms on the unit interval are exactly models of
piecewise polynomial functions with degrees not larger than $0$. In \cite%
{saum:12}, it is shown that models of piecewise polynomial functions have
nice analytical and statistical properties. Let us recall two of them.

In Lemma 8 of \cite{saum:12}, it is proved that if the distribution $P^{X}$
has a density with respect to the Lebesgue measure $%
\leb%
$ on $\mathcal{X=}\left[ 0,1\right] $ which is uniformly bounded away from
zero and if the considered partition $\mathcal{P}$ is lower regular with
respect to $%
\leb%
$ - that is there exists a positive constant $c$ such that $\left\vert 
\mathcal{P}\right\vert \inf_{I\in \mathcal{P}}%
\leb%
\left( I\right) \geq c>0$ - then the associated model of piecewise
polynomial functions is equipped with a localized orthonormal basis in $%
L_{2}\left( P^{X}\right) $. For a formal definition of a localized basis,
see Section \ref{section_proof_slope_reg}\ below. Since the pioneering work
of Birg\'{e} and Massart 
\cite{BirgeMassart:93,BirgeMassart:98,Massart:07}%
, the property of localized basis is known to play a key role in
M-estimation and model selection using vector spaces or more general sieves.

Considering models of piecewise polynomial functions on the unit interval,
where the density of $P^{X}$ with respect to $%
\leb%
$ is both uniformly bounded and bounded away from 0 and where the underlying
partition is lower regular with respect to $%
\leb%
$, it is shown in Lemma 9 of \cite{saum:12} that the least-squares estimator 
$s_{n}\left( M\right) $ converges in sup-norm to the linear projection $%
s_{M} $ of the regression function $s_{\ast }$.

Assumptions of lower regularity of the considered partitions as well as the
existence of a uniformly bounded density of $P^{X}$ with respect to the
Lebesgue measure on $\mathcal{X}$, will thus naturally arise when dealing
with least-squares model selection using piecewise polynomial functions -
see Section \ref{section_main_assumptions} below. Furthermore, the
interested reader will find in Section \ref{section_proof_slope_reg} a more
general version of our results, available for linear models equipped with a
localized basis and where least-squares estimators converge in sup-norm to
the linear projections of the regression function onto the models.

\section{The slope heuristics\label{section_slope_heuristics}}

\subsection{Underlying concepts\label{section_underlying_concepts}}

In order to clarify our approach and to highlight the connection of the
present paper with the results previously established in \cite{saum:12}, we
first give a brief heuristic explanation of the major mathematical facts
underlying the slope phenomenon.

We rewrite the definition of the oracle model $M_{\ast }$ given in (\ref%
{oracle_model}). For any $M\in \mathcal{M}_{n}$, the excess loss $\ell
\left( s_{\ast },s_{n}\left( M\right) \right) =P\left( Ks_{n}\left( M\right)
\right) -P\left( Ks_{\ast }\right) $ is the difference between the loss of
the estimator $s_{n}\left( M\right) $ and the loss of the target $s_{\ast }$%
. As $P\left( Ks_{\ast }\right) $ is independent of $M$ varying in $\mathcal{%
M}_{n}$, it holds%
\begin{align*}
M_{\ast }& \in \arg \min_{M\in \mathcal{M}_{n}}\left\{ P\left( Ks_{n}\left(
M\right) \right) \right\} \\
& =\arg \min_{M\in \mathcal{M}_{n}}\left\{ P_{n}\left( Ks_{n}\left( M\right)
\right) +%
\pen%
_{\text{id}}\left( M\right) \right\} \text{ ,}
\end{align*}%
where for all $M\in \mathcal{M}_{n}$, 
\begin{equation*}
\pen%
_{\text{id}}\left( M\right) :=P\left( Ks_{n}\left( M\right) \right)
-P_{n}\left( Ks_{n}\left( M\right) \right) \text{ .}
\end{equation*}%
The penalty function $%
\pen%
_{\text{id}}$ is called the \textit{ideal penalty} - as it allows to select
the oracle - and is unknown because it depends on the distribution of the
data. As pointed out by Arlot and Massart \cite{ArlotMassart:09}, the main
idea of penalization in the efficiency problem is to give some sharp
estimate, up to a constant, of the ideal penalty. This would yield an
(asymptotically) unbiased - or uniformly biased over the collection of
models $\mathcal{M}_{n}$ - estimation of the loss. Such a penalization would
lead to a sharp oracle inequality for the selected model.

A penalty term $%
\pen%
_{\text{opt}}$ is said to be optimal if it achieves an oracle inequality
with leading constant converging to one when the sample size $n$ tends to
infinity.

Concerning the estimation of the optimal penalty, Arlot and Massart \cite%
{ArlotMassart:09} conjectured that the mean of the empirical excess loss $%
\mathbb{E}\left[ P_{n}\left( Ks_{M}-Ks_{n}\left( M\right) \right) \right] $
satisfies the following slope heuristics in a quite general M-estimation
framework:

\begin{description}
\item[(\textbf{i})] If a penalty $%
\pen%
:\mathcal{M}_{n}\longrightarrow \mathbb{R}_{+}$ is such that, for all models 
$M\in \mathcal{M}_{n}$, 
\begin{equation*}
\pen%
\left( M\right) \leq \left( 1-\delta \right) \mathbb{E}\left[ P_{n}\left(
Ks_{M}-Ks_{n}\left( M\right) \right) \right]
\end{equation*}%
with $\delta >0$, then the dimension of the selected model $%
\widehat{M}%
$ is \textquotedblleft very large\textquotedblright\ and the excess loss of
the selected estimator $s_{n}\left( 
\widehat{M}%
\right) $ is \textquotedblleft much larger\textquotedblright\ than the
excess loss of the oracle.

\item[(\textbf{ii})] If $%
\pen%
\approx \left( 1+\delta \right) \mathbb{E}\left[ P_{n}\left(
Ks_{M}-Ks_{n}\left( M\right) \right) \right] $ with $\delta >0$, then the
corresponding model selection procedure satisfies an oracle inequality with
a leading constant $C\left( \delta \right) <+\infty $ and the dimension of
the selected model is \textquotedblleft not too large\textquotedblright .
Moreover, 
\begin{equation*}
\pen%
_{\text{opt}}\left( M\right) \approx 2\mathbb{E}\left[ P_{n}\left(
Ks_{M}-Ks_{n}\left( M\right) \right) \right]
\end{equation*}%
is an optimal penalty.
\end{description}

\noindent The mean of the empirical excess loss on $M$, when $M$ varies in $%
\mathcal{M}_{n}$, is thus conjectured to be the maximal value of penalty
under which the model selection procedure totally misbehaves or,
equivalently, the minimum value of penalty above which the procedure
achieves an oracle inequality. It is called the \textit{minimal penalty, }%
denoted by $%
\pen%
_{\min }$:%
\begin{equation*}
\text{for all }M\in \mathcal{M}_{n}\text{, \ \ \ \ \ }%
\pen%
_{\min }\left( M\right) =\mathbb{E}\left[ P_{n}\left( Ks_{M}-Ks_{n}\left(
M\right) \right) \right] \text{ .}
\end{equation*}%
The optimal penalty is then close to twice the minimal one,%
\begin{equation}
\pen%
_{\text{opt}}\approx 2%
\pen%
_{\min }\text{ .}  \label{relation_pen_opt_pen_min}
\end{equation}%
Let us now briefly explain the points (i) and (ii) above. We give in Section %
\ref{section_main_theorems_reg} precise results which validate the slope
heuristics for models of piecewise polynomial functions.

If the chosen penalty is less than the minimal one, $%
\pen%
=\left( 1-\delta \right) 
\pen%
_{\min }$ with $\delta \in \left[ 0,1\right] $, the algorithm minimizes over 
$\mathcal{M}_{n}$,%
\begin{align*}
& P_{n}\left( Ks_{n}\left( M\right) \right) +%
\pen%
\left( M\right) -P_{n}\left( Ks_{\ast }\right) \\
& =P\left( Ks_{M}-Ks_{\ast }\right) +\left( P_{n}-P\right) \left(
Ks_{M}-Ks_{\ast }\right) -P_{n}\left( Ks_{M}-Ks_{n}\left( M\right) \right) +%
\pen%
\left( M\right) \\
& =P\left( Ks_{M}-Ks_{\ast }\right) +\left( P_{n}-P\right) \left(
Ks_{M}-Ks_{\ast }\right) -\delta P_{n}\left( Ks_{M}-Ks_{n}\left( M\right)
\right) \\
& \text{ \ \ \ }+\left( 1-\delta \right) \left( \mathbb{E}\left[ P_{n}\left(
Ks_{M}-Ks_{n}\left( M\right) \right) \right] -P_{n}\left(
Ks_{M}-Ks_{n}\left( M\right) \right) \right) \\
& \approx \ell \left( s_{\ast },s_{M}\right) -\delta P_{n}\left(
Ks_{M}-Ks_{n}\left( M\right) \right) \text{ .}
\end{align*}%
In the latter identity, we neglect the difference between the empirical and
true loss of the projections $s_{M}$ and the deviations of the empirical
excess loss $P_{n}\left( Ks_{M}-Ks_{n}\left( M\right) \right) $. Indeed, as
shown by Boucheron and Massart \cite{BouMas:10}, the empirical excess loss
satisfies a concentration inequality in a general framework, which allows to
neglect the difference with its mean, at least for models that are not too
small.

As the empirical excess loss is increasing and the excess loss of the
projection $s_{M}$ is decreasing with respect to the complexity of the
models, the penalized criterion is (almost) decreasing with respect to the
complexity of the models, and the selected model is among the largest of the
collection.

On the contrary, if the chosen penalty is greater than the minimal one, $%
\pen%
=\left( 1+\delta \right) 
\pen%
_{\min }$ with $\delta >0$, then by the same kind of manipulations, the
selected model minimizes the following criterion, for all $M\in \mathcal{M}%
_{n}$,%
\begin{equation}
P_{n}\left( Ks_{n}\left( M\right) \right) +%
\pen%
\left( M\right) -P_{n}\left( Ks_{\ast }\right) \approx \ell \left( s_{\ast
},s_{M}\right) +\delta P_{n}\left( Ks_{M}-Ks_{n}\left( M\right) \right) 
\text{ .}  \label{criterion_heur}
\end{equation}%
The selected model thus achieves a trade-off between the bias of the models
which decreases with the complexity and the empirical excess loss which
increases with the complexity of the models. The selected dimension would
then be reasonable, and the trade-off between the bias and the complexity of
the models is likely to give some oracle inequality.

Finally, if we take $\delta =1$ in the latter case, $%
\pen%
=2\times 
\pen%
_{\min }$, and if we assume that the empirical excess loss is equivalent to
the excess loss,%
\begin{equation}
P_{n}\left( Ks_{M}-Ks_{n}\left( M\right) \right) \sim P\left( Ks_{n}\left(
M\right) -Ks_{M}\right) \text{ ,}  \label{equivalence_p1_p2}
\end{equation}%
then according to (\ref{criterion_heur}) the selected model almost minimizes%
\begin{equation*}
P\left( Ks_{M}-Ks_{\ast }\right) +P_{n}\left( Ks_{M}-Ks_{n}\left( M\right)
\right) \approx \ell \left( s_{\ast },s_{M}\right) +P\left( Ks_{n}\left(
M\right) -Ks_{M}\right) \approx \ell \left( s_{\ast },s_{n}\left( M\right)
\right) \text{ .}
\end{equation*}%
Hence,%
\begin{equation*}
\ell \left( s_{\ast },s_{n}\left( 
\widehat{M}%
\right) \right) \approx \ell \left( s_{\ast },s_{n}\left( M_{\ast }\right)
\right)
\end{equation*}%
and the procedure is nearly optimal.

One can find in \cite{saum:12} some results showing that (\ref%
{equivalence_p1_p2}) is a quite general fact in least-squares regression and
is in particular satisfied when considering models of piecewise polynomial
functions. Thus, these results represent a preliminary material for the
present study, and we shall base our arguments on the results exposed in 
\cite{saum:12}.

\subsection{Assumptions and comments\label{section_main_assumptions}}

\noindent We take $\mathcal{X=}\left[ 0,1\right] $, $%
\leb%
$ is the Lebesgue measure on $\mathcal{X}$, and linear models $M\in \mathcal{%
M}_{n}$ are models of piecewise polynomial functions. We denote by $\mathcal{%
P}_{M}$ the partition of $\mathcal{X}$ underlying the model $M$.

\bigskip

\noindent \textbf{Set of assumptions for piecewise polynomial functions}: (%
\textbf{SAPP})

\bigskip

\begin{description}
\item[(\textbf{P1})] there exist two positive constants $c_{\mathcal{M}%
},\alpha _{\mathcal{M}}$ such that $%
\card%
\left( \mathcal{M}_{n}\right) \leq c_{\mathcal{M}}n^{\alpha _{\mathcal{M}}}$ 
$.$

\item[(\textbf{P2})] there exists a positive constant $A_{\mathcal{M},+}$
such that for every $M\in \mathcal{M}_{n},$ $1\leq D_{M}\leq A_{\mathcal{M}%
,+}n\left( \ln n\right) ^{-2}\leq n$ $.$

\item[(\textbf{P3})] there exist $c_{rich}>0$, $A_{rich}>0$ and $%
M_{0},M_{1}\in \mathcal{M}_{n}$ such that $D_{M_{0}}\in \left[ n^{1/\left(
1+\beta _{+}\right) },c_{rich}n^{1/\left( 1+\beta _{+}\right) }\right] $ and 
$D_{M_{1}}\geq A_{rich}n\left( \ln n\right) ^{-2}$, where $\beta _{+}$ is
defined in (\textbf{Ap}$_{u}$).

\item[(\textbf{Ap}$_{u}$)] there exist $\beta _{+}>0$ and $C_{+}>0$ such that%
\begin{equation*}
\ell \left( s_{\ast },s_{M}\right) \leq C_{+}D_{M}^{-\beta _{+}}\text{ }.
\end{equation*}

\item[(\textbf{An})] There exists a constant $\sigma _{\min }$ such that $%
\sigma \left( X_{i}\right) \geq \sigma _{\min }>0$ $a.s.$

\item[(\textbf{Ab})] There exists a positive constant $A$, that bounds the
data: $\left\vert Y_{i}\right\vert \leq A<\infty .$

\item[(\textbf{Ad}$_{%
\leb%
}$)] $P^{X}$ has a density $f$ with respect to $%
\leb%
$ satisfying for some constants $c_{\min }$ and $c_{\max }$, that%
\begin{equation*}
\text{\ }0<c_{\min }\leq f\left( x\right) \leq c_{\max }<\infty ,\text{ \ }%
\forall x\in \left[ 0,1\right] \text{ }.
\end{equation*}

\item[(\textbf{Aud})] there exists $r\in \mathbb{N}^{\ast }$ such that, for
all $M\in \mathcal{M}_{n}$, all $I\in \mathcal{P}_{M}$ and all $p\in M$,%
\begin{equation*}
\deg \left( p_{\mid I}\right) \leq r\text{ }.
\end{equation*}

\item[(\textbf{Alr})] a positive constant $c_{\mathcal{M},%
\leb%
}$ exists such that, for all $M\in \mathcal{M}_{n}$,%
\begin{equation*}
0<c_{\mathcal{M},%
\leb%
}\leq \left\vert \mathcal{P}_{M}\right\vert \inf_{I\in \mathcal{P}_{M}}%
\leb%
\left( I\right) <+\infty \text{ }.
\end{equation*}
\end{description}

\bigskip

\noindent The set of assumptions (\textbf{SAPP}) can be divided into three
groups. Firstly, assumptions (\textbf{P1}), (\textbf{P2}), (\textbf{P3}) and
(\textbf{Ap}$_{u}$) are linked to properties of the collection of models $%
\mathcal{M}_{n}$. Secondly, assumptions (\textbf{An}), (\textbf{Ab}) and (%
\textbf{Ad}$_{%
\leb%
}$) give some constraints on the general regression relation stated in (\ref%
{regression_model_2}). Thirdly, assumptions (\textbf{Aud}) and (\textbf{Alr}%
) specify some quantities related to the choice of the models of piecewise
polynomial functions.

Assumption (\textbf{P1}) states that the collection of models has a
\textquotedblleft small\textquotedblright\ complexity, more precisely a
polynomially increasing one with respect to the amount of data. For this
kind of complexities, if one wants to design a good model selection
procedure for prediction, the chosen penalty should estimate the mean of the
ideal one on each model, up to a constant. Indeed, as Talagrand's type
concentration inequalities for the empirical process are exponential, they
allow to neglect the deviations of the quantities of interest from their
mean, uniformly over the collection of models. This is not the case for
large collections of models, where one has to put an extra-log factor
depending on the complexity of the collection of models inside the penalty,
see for instance 
\cite{BirMassart:97,BarBirMassart:99}%
.

We assume in (\textbf{P3}) that the collection of models contains a model $%
M_{0}$ of reasonably large dimension and a model $M_{1}$ of high dimension,
which is necessary since we prove the existence of a jump between high and
reasonably large dimensions. One can notice that in practice, the parameter $%
\beta _{+}$, which depends on the bias of the model is not known and so the
existence of $M_{0}$ is not straightforward. However, it suffices for the
statistician to take at least one model per dimension lower than the chosen
upper bound to ensure the existence of $M_{0}$ and $M_{1}$.

We require in (\textbf{Ap}$_{u}$) for the quality of approximation of the
collection of models to be good enough in terms of the quadratic loss. More
precisely, we ask for a polynomial decrease of excess loss of linear
projections of the regression function onto the models. It is well-known
that piecewise polynomial functions uniformly bounded in their degrees have
good approximation properties in Besov spaces. More precisely, as stated in
Lemma 12 of Barron, Birg\'{e} and Massart \cite{BarBirMassart:99}, if $%
\mathcal{X=}\left[ 0,1\right] $ and the regression function $s_{\ast }$
belongs to the Besov space $B_{\alpha ,p,\infty }\left( \mathcal{X}\right) $
(see the definition in \cite{BarBirMassart:99}), then taking models of
piecewise polynomial functions of degree bounded by $r>\alpha -1$ on regular
partitions with respect to the Lebesgue measure $%
\leb%
$ on $\mathcal{X}$, and assuming that $P^{X}$ has a density with respect to $%
\leb%
$ which is bounded in sup-norm, assumption (\textbf{Ap}$_{u}$) is satisfied.

Assumption (\textbf{Ab}) is rather restrictive, since it excludes Gaussian
noise. However, the assumption of bounded noise is somehow classical when
dealing with M-estimation and related procedures. Indeed, a central tool in
this field is empirical process theory and more especially, concentration
inequalities for the supremum of the empirical process. We used the
classical inequalities of Bousquet, and Klein and Rio in \cite{saum:12}. As
a matter of fact, we do not know yet if an adaptation of our proofs
(including results established in \cite{saum:12}) by using extensions of the
latter inequalities to some unbounded cases - as for instance in Adamczak's
concentration inequalities \cite{Adamczak:08} - would be possible.

The noise restriction stated in (\textbf{An}) is needed to derive our
results which are optimal to the first order. More precisely, it allows in 
\cite{saum:12} to obtain sharp lower bounds for the true and empirical
excess losses on a fixed model. This assumption is also needed in the work
of Arlot and Massart \cite{ArlotMassart:09} concerning the case of histogram
models. As it is noticed in Section 5.3 of \cite{saum:12}, assumption (%
\textbf{An}) could be replaced by the following assumption, which states
that the partitions underlying the models of piecewise polynomial functions
are regular from above with respect to the Lebesgue measure on $\left[ 0,1%
\right] $.

\begin{description}
\item[(Aur)] a positive constant $c_{\mathcal{M},%
\leb%
}^{+}$ exists such that, for all $M\in \mathcal{M}_{n}$,%
\begin{equation*}
\left\vert \mathcal{P}_{M}\right\vert \sup_{I\in \mathcal{P}_{M}}%
\leb%
\left( I\right) \leq c_{\mathcal{M},%
\leb%
}^{+}\text{ .}
\end{equation*}
\end{description}

Assumptions (\textbf{Ad}$_{%
\leb%
}$), (\textbf{Aud}) and (\textbf{Alr}) imply several important properties
for the models of piecewise polynomial functions, such as the existence of
an orthonormal localized basis in each model or the consistency in sup-norm
of least-squares estimators toward the projections of the target onto the
models. See also Sections \ref{section_piecewise_polynomials}\ and \ref%
{section_more_general}\ for further comments about these properties.

\subsection{Statement of the theorems\label{section_main_theorems_reg}}

We are now able to state our main results leading to the slope heuristics.
They describe the behavior of the penalization procedure defined in (\ref%
{def_proc_2_reg}).

\begin{theorem}
\label{theorem_min_pen_reg_pp}Take a positive penalty: for all $M\in 
\mathcal{M}_{n}$, $%
\pen%
\left( M\right) \geq 0$. Suppose that the assumptions (\textbf{SAPP}) of
Section \ref{section_main_assumptions} hold, and furthermore suppose that
for $A_{%
\pen%
}\in \left[ 0,1\right) $ and $A_{p}>0$ the model $M_{1}$ of assumption (%
\textbf{P3}) satisfies 
\begin{equation}
\text{ }0\leq 
\pen%
\left( M_{1}\right) \leq A_{%
\pen%
}\mathbb{E}\left[ P_{n}\left( Ks_{M_{1}}-Ks_{n}\left( M_{1}\right) \right) %
\right] \text{ },  \label{majo_pen_pp}
\end{equation}%
with probability at least $1-A_{p}n^{-2}$. Then there exist a constant $%
A_{1}>0$ only depending on constants in \textit{(\textbf{SAPP})}, as well as
an integer $n_{0}$ and a positive constant $A_{2}$ only depending on $A_{%
\pen%
}$ and on constants in \textit{(\textbf{SAPP})} such that, for all $n\geq
n_{0}$, it holds with probability at least $1-A_{1}n^{-2}$,%
\begin{equation*}
D_{%
\widehat{M}%
}\geq A_{2}n\ln \left( n\right) ^{-2}
\end{equation*}%
and%
\begin{equation}
\ell \left( s_{\ast },s_{n}\left( 
\widehat{M}%
\right) \right) \geq \frac{n^{\beta _{+}/\left( 1+\beta _{+}\right) }}{%
\left( \ln n\right) ^{3}}\inf_{M\in \mathcal{M}_{n}}\left\{ \ell \left(
s_{\ast },s_{n}\left( M\right) \right) \right\} \text{ },
\label{bad_oracle_min_pen_pp}
\end{equation}%
where $\beta _{+}>0$ is defined in assumption (\textbf{Ap}$_{u}$) of (%
\textbf{SAPP}).
\end{theorem}

\noindent Theorem \ref{theorem_min_pen_reg_pp} justifies the first part (%
\textbf{i}) of the slope heuristics exposed in Section \ref%
{section_slope_heuristics}. As a matter of fact, it shows that there exists
a level such that, if the penalty is smaller than this level for one of the
largest models, then the dimension of the output is among the largest
dimensions of the collection and the excess loss of the selected estimator
is much larger than the excess loss of the oracle. Moreover, this level is
given by the mean of the empirical excess loss of the least-squares
estimator on each model. Let us also notice that the lower bound given in (%
\ref{bad_oracle_min_pen_pp}) gets worse as $\beta _{+}$ increases. This is
due to the fact that when $\beta _{+}$ increases, the approximation
properties of the models improve and the performances in terms of excess
loss for the oracle estimator also improve.

\noindent The following theorem validates the second part of the slope
heuristics.

\begin{theorem}
\label{theorem_opt_pen_reg_pp}Suppose that the assumptions \textit{(\textbf{%
SAPP})} of Section \ref{section_main_assumptions} hold, and furthermore
suppose that for some $\delta \in \left[ 0,1\right) $ and $A_{p},A_{r}>0$,
there exists an event of probability at least $1-A_{p}n^{-2}$ on which, for
every model $M\in \mathcal{M}_{n}$ such that $D_{M}\geq A_{\mathcal{M}%
,+}\left( \ln n\right) ^{3}$, it holds%
\begin{equation}
\left\vert 
\pen%
\left( M\right) -2\mathbb{E}\left[ P_{n}\left( Ks_{M}-Ks_{n}\left( M\right)
\right) \right] \right\vert \leq \delta \left( \ell \left( s_{\ast
},s_{M}\right) +\mathbb{E}\left[ P_{n}\left( Ks_{M}-Ks_{n}\left( M\right)
\right) \right] \right)  \label{pen_id_pp}
\end{equation}%
together with%
\begin{equation}
\left\vert 
\pen%
\left( M\right) \right\vert \leq A_{r}\left( \frac{\ell \left( s_{\ast
},s_{M}\right) }{\left( \ln n\right) ^{2}}+\frac{\left( \ln n\right) ^{3}}{n}%
\right) \text{ .}  \label{pen_id_2_pp}
\end{equation}%
Then, for any $\eta \in \left( 0,\beta _{+}/\left( 1+\beta _{+}\right)
\right) $, there exist an integer $n_{0}$ only depending on $\eta ,\delta $
and $\beta _{+}$ and on constants in \textit{(\textbf{SAPP})}, a positive
constant $A_{3}$ only depending on $c_{\mathcal{M}}$ given in \textit{(%
\textbf{SAPP})} and on $A_{p}$, two positive constants $A_{4}$ and $A_{5}$
only depending on constants in \textit{(\textbf{SAPP})} and on $A_{r}$ and a
sequence 
\begin{equation}
\theta _{n}\leq \frac{A_{4}}{\left( \ln n\right) ^{1/4}}
\label{def_theta_n_pp}
\end{equation}%
such that it holds for all $n\geq n_{0}$, with probability at least $%
1-A_{3}n^{-2}$, 
\begin{equation*}
D_{%
\widehat{M}%
}\leq n^{\eta +1/\left( 1+\beta _{+}\right) }
\end{equation*}%
and%
\begin{equation}
\ell \left( s_{\ast },s_{n}\left( 
\widehat{M}%
\right) \right) \leq \left( \frac{1+\delta }{1-\delta }+\frac{5\theta _{n}}{%
\left( 1-\delta \right) ^{2}}\right) \ell \left( s_{\ast },s_{n}\left(
M_{\ast }\right) \right) +A_{5}\frac{\left( \ln n\right) ^{3}}{n}\text{ .}
\label{oracle_opt_gene_pp}
\end{equation}%
Assume that in addition, the following assumption holds,

\begin{description}
\item[(\textbf{Ap})] The bias decreases like a power of $D_{M}$: there exist 
$\beta _{-}\geq \beta _{+}>0$ and $C_{+},C_{-}>0$ such that%
\begin{equation*}
C_{-}D_{M}^{-\beta _{-}}\leq \ell \left( s_{\ast },s_{M}\right) \leq
C_{+}D_{M}^{-\beta _{+}}\text{ }.
\end{equation*}
\end{description}

\noindent Then it holds for all $n\geq n_{0}\left( \left( \text{\textbf{SAPP}%
}\right) ,C_{-},\beta _{-},\beta _{+},\eta ,\delta \right) $, with
probability at least $1-A_{3}n^{-2}$,%
\begin{equation}
A_{\mathcal{M},+}\left( \ln n\right) ^{3}\leq D_{%
\widehat{M}%
}\leq n^{\eta +1/\left( 1+\beta _{+}\right) }  \label{dim_ap_pp}
\end{equation}%
and%
\begin{equation}
\ell \left( s_{\ast },s_{n}\left( 
\widehat{M}%
\right) \right) \leq \left( \frac{1+\delta }{1-\delta }+\frac{5\theta _{n}}{%
\left( 1-\delta \right) ^{2}}\right) \ell \left( s_{\ast },s_{n}\left(
M_{\ast }\right) \right) \text{ .}  \label{oracle_opt_pp}
\end{equation}
\end{theorem}

Theorem \ref{theorem_opt_pen_reg_pp} states that if the penalty is close to
twice the minimal one, then the selected estimator satisfies a pathwise
oracle inequality with constant almost one, and so the model selection
procedure is approximately optimal. Moreover, the dimension of the selected
model is of reasonable dimension, bounded by a power less than one of the
sample size.

Condition (\textbf{Ap}) allows to remove the remainder terms from the oracle
inequality (\ref{oracle_opt_gene_pp}) by ensuring that the selected model is
of dimension not too small, as stated in (\ref{dim_ap_pp}). Assumption (%
\textbf{Ap}) is the conjunction of assumption (\textbf{Ap}$_{u}$) with a
polynomial lower bound of the bias of the models. On histogram models, Arlot
showed in Section 8.10 of \cite{Arlot:07} that this lower bound is satisfied
for non constant $\alpha $-H\"{o}lder, $\alpha \in \left( 0,1\right] $,
regression functions and for regular partitions.

Finally, from Theorems \ref{theorem_min_pen_reg_pp} and \ref%
{theorem_opt_pen_reg_pp}, we identify the minimal penalty with the mean of
the empirical excess loss on each model, 
\begin{equation*}
\pen%
_{\min }\left( M\right) =\mathbb{E}\left[ P_{n}\left( Ks_{M}-Ks_{n}\left(
M\right) \right) \right] \text{ ,}
\end{equation*}%
thus generalizing the results of Arlot and Massart in \cite{ArlotMassart:09}
to the case of piecewise polynomial functions.

\section{Hold-out penalization\label{section_hold_out}}

The conditions on the penalty given in Theorems \ref{theorem_min_pen_reg_pp}
and \ref{theorem_opt_pen_reg_pp} can not be directly checked in practice.
Indeed, they are expressed in terms of the mean of the empirical excess loss
on each model, which is an unknown quantity in general. Nevertheless, in the
homoscedastic case, it is easy to see that Mallows' penalty is a
nonasymptotic quasi-optimal penalty. According to Theorem \ref%
{theorem_opt_pen_reg_pp}, such a penalty is given by twice the mean of the
empirical excess loss. Now, using Theorem 10 of \cite{saum:12}, we get (with
an explicit control of the second order terms in the following equivalence),%
\begin{equation*}
2\mathbb{E}\left[ P_{n}\left( Ks_{M}-Ks_{n}\left( M\right) \right) \right]
\sim \frac{1}{2}\mathcal{K}_{1,M}^{2}\frac{D_{M}}{n}\text{ ,}
\end{equation*}%
where $\mathcal{K}_{1,M}^{2}=1/D_{M}\sum_{k=1}^{D_{M}}\mathbb{E}\left(
\left( \psi _{1,M}\left( X,Y\right) \cdot \varphi _{k}\left( X\right)
\right) ^{2}\right) $, $\psi _{1,M}\left( X,Y\right) =-2\left( Y-s_{M}\left(
X\right) \right) ~$and $\left( \varphi _{k}\right) _{k=1}^{D_{M}}$ is an
orthonormal basis in $\left( M,\left\Vert \cdot \right\Vert _{2}\right) $.
By easy computations, we deduce that if the noise is homoscedastic, that is $%
\sigma ^{2}\left( X\right) \equiv \sigma ^{2}>0$, it holds%
\begin{equation}
\frac{1}{2}\mathcal{K}_{1,M}^{2}\frac{D_{M}}{n}=2\sigma ^{2}\frac{D_{M}}{n}+%
\mathbb{E}\left[ \left( s_{\ast }-s_{M}\right) ^{2}\frac{\sum_{i=1}^{D_{M}}%
\varphi _{k}^{2}}{n}\right] \text{ .}  \label{homo}
\end{equation}%
The second term at the right of identity (\ref{homo}) being negligible for
models of interest in the conditions of Theorem \ref{theorem_opt_pen_reg_pp}
(thanks to Lemma 7 in \cite{saum:12}, which implies that $%
\sum_{i=1}^{D_{M}}\varphi _{k}^{2}\leq LD_{M}$ for some constant $L>0$), we
conclude that an\ asymptotically optimal penalty is given by $2\sigma
^{2}D_{M}/n$, which is Mallows' classical penalty.

In the case where the noise level is homoscedastic but unknown, Mallows'
penalty is only known through a constant, the noise level, which can be
estimated \textit{via} the slope heuristics (for practical issues about the
slope heuristics, see Baudry et al. \cite{BauMauMich:12}). But in the common
situation where the noise level is sufficiently heteroscedastic, the shape
of the ideal penalty is not linear in the dimension of the models and not
even a \textit{function }of the linear dimensions. In such a case, Arlot 
\cite{Arl:2010} proved that any calibration of a linear penalty leads to a
suboptimal procedure, but yet can achieve an oracle inequality with a
leading constant more than one.

In order to achieve a nearly optimal selection procedure in the general
situation, it remains to estimate the ideal penalty or, thanks to the slope
heuristics, the shape of the ideal penalty. This section is devoted to this
task. We propose a hold-out type penalty that automatically adapts to
heteroscedasticity.\ Let us now detail our hold-out penalization procedure.

The ideal penalty is defined by%
\begin{equation*}
\pen%
_{\text{id}}\left( M\right) :=P\left( Ks_{n}\left( M\right) \right)
-P_{n}\left( Ks_{n}\left( M\right) \right) \text{ ,}
\end{equation*}%
for all $M\in \mathcal{M}_{n}$. A natural idea is to divide the data into
two groups, indexed by $I_{1}$ and $I_{2}$, satisfying $I_{1}\cap
I_{2}=\emptyset $ and $I_{1}\cup I_{2}=\left\{ 1,...,n\right\} $ and to
propose the following hold-out type penalty,%
\begin{equation*}
\pen%
_{ho,C}\left( M\right) :=C\left( P_{n_{2}}\left( Ks_{n_{1}}\left( M\right)
\right) -P_{n_{1}}\left( Ks_{n_{1}}\left( M\right) \right) \right) \text{ ,}
\end{equation*}%
where $P_{n_{i}}=1/n_{i}\sum_{j\in I_{i}}\delta _{\xi _{j}}$, $n_{i}=$Card$%
\left( I_{i}\right) $, for $i=1,2$, $s_{n_{1}}\left( M\right) \in \arg
\min_{s\in M}P_{n_{1}}\left( Ks\right) $ and $C>0$ is a constant to be
determined. Indeed, if $n_{1}$ is not too small, $P_{n_{1}}\left(
Ks_{n_{1}}\left( M\right) \right) $ is likely to vary like $P_{n}\left(
Ks_{n}\left( M\right) \right) $ and $P_{n_{2}}\left( Ks_{n_{1}}\left(
M\right) \right) $ is, conditionally to $\left( \xi _{j}\right) _{j\in
I_{1}} $, an unbiased estimate of $P\left( Ks_{n_{1}}\left( M\right) \right) 
$, which again is likely to vary like $P\left( Ks_{n}\left( M\right) \right) 
$. Moreover, we see from Theorem 10 in \cite{saum:12} that when the model $M$
is fixed, the quantities $P_{n}\left( Ks_{n}\left( M\right) \right) $ and $%
P\left( Ks_{n}\left( M\right) \right) $ are almost inversely proportional to 
$n$, so a good constant in front of the hold-out penalty should be $%
C_{opt}=n_{1}/n$.

The previous observation is justified by the following theorem, where for
the sake of clarity we fixed $n_{1}=n_{2}=n/2$. For a more general version
of Theorem \ref{theorem_pen_n1n2_pp}, see Section \ref%
{section_proof_hold_out}. We set 
\begin{equation}
\pen%
_{ho}\left( M\right) =\frac{1}{2}\left( P_{n_{2}}\left( Ks_{n_{1}}\left(
M\right) \right) -P_{n_{1}}\left( Ks_{n_{1}}\left( M\right) \right) \right) 
\text{ \ and \ }\widehat{M}_{1/2}\in \arg \min_{M\in \mathcal{M}_{n}}\left\{
P_{n}\left( Ks_{n}\left( M\right) \right) +%
\pen%
_{ho}\left( M\right) \right\} \text{ .}  \label{def_proc_half_and_half}
\end{equation}

\begin{theorem}
\label{theorem_pen_n1n2_pp}Consider the procedure defined in (\ref%
{def_proc_half_and_half}), with $n_{1}=n_{2}=n/2$. Suppose that the
assumptions \textit{(\textbf{SAPP})} of Section \ref%
{section_main_assumptions} hold. Then, for any $\eta \in \left( 0,\beta
_{+}/\left( 1+\beta _{+}\right) \right) $, there exist an integer $n_{0}$
only depending on $\eta $ and on constants in \textit{(\textbf{SAPP})}, a
positive constant $A_{6}$ only depending on $c_{\mathcal{M}}$ given in 
\textit{(\textbf{SAPP})}, two positive constants $A_{7}$ and $A_{8}$ only
depending on constants in \textit{(\textbf{SAPP})} and a sequence $\theta
_{n}\leq A_{7}\left( \ln n\right) ^{-1/4}$ such that it holds for all $n\geq
n_{0}$, with probability at least $1-A_{6}n^{-2}$, 
\begin{equation*}
D_{%
\widehat{M}%
_{1/2}}\leq n^{\eta +1/\left( 1+\beta _{+}\right) }
\end{equation*}%
and 
\begin{equation}
\ell \left( s_{\ast },s_{n}\left( 
\widehat{M}%
_{1/2}\right) \right) \leq \left( 1+\theta _{n}\right) \ell \left( s_{\ast
},s_{n}\left( M_{\ast }\right) \right) +A_{8}\frac{\left( \ln n\right) ^{3}}{%
n}\text{ .}  \label{oracle_opt_gene_n1_pp}
\end{equation}%
Assume that in addition (\textbf{Ap}) holds (see Theorem \ref%
{theorem_opt_pen_reg_pp}). Then it holds for all $n\geq n_{0}\left( \left( 
\text{\textbf{SAPP}}\right) ,C_{-},\beta _{-},\eta \right) $, with
probability at least $1-A_{6}n^{-2}$,%
\begin{equation*}
A_{\mathcal{M},+}\left( \ln n\right) ^{3}\leq D_{%
\widehat{M}%
_{1/2}}\leq n^{\eta +1/\left( 1+\beta _{+}\right) }
\end{equation*}%
and%
\begin{equation}
\ell \left( s_{\ast },s_{n}\left( 
\widehat{M}%
_{1/2}\right) \right) \leq \left( 1+\theta _{n}\right) \inf_{M\in \mathcal{M}%
_{n}}\left\{ \ell \left( s_{\ast },s_{n}\left( M\right) \right) \right\} 
\text{ .}  \label{oracle_opt_n1_pp}
\end{equation}
\end{theorem}

\bigskip

\noindent Theorem \ref{theorem_pen_n1n2_pp} shows the asymptotic optimality
of the hold-out penalization procedure, for a half-and-half split of the
data. This is a remarkable fact compared to the classical hold-out, defined
by%
\begin{equation}
\widehat{M}_{ho}\in \arg \min_{M\in \mathcal{M}_{n}}\left\{ P_{n_{2}}\left(
Ks_{n_{1}}\left( M\right) \right) \right\} \text{ .}  \label{def_ho}
\end{equation}%
Indeed, the choice $n_{1}=n/2$ in (\ref{def_ho}) is likely to lead to an
asymptotically suboptimal procedure, as the criterion is close in
expectation to $P\left( Ks_{n/2}\left( M\right) \right) $, and so is close
to the oracle, but for $n/2$ data points. The hold-out penalization allows
us to overcome this difficulty. Arlot 
\cite{Arl:2008a,Arl:09}
described similar advantages for resampling and $V$-fold penalties.

Notice also that the random hold-out penalty proposed by Arlot \cite{Arl:09}
is proportional to the mean along the splits of our hold-out penalty,
providing thus a \textquotedblleft stabilization effect\textquotedblright\
in practice. This should bring some improvement compared to our unique
split, at the price of increased computational cost. However, the
stabilization effect seems more difficult to study mathematically, and our
results provide a first step toward the study of the more complicated
resampling penalties.

\section{Proofs\label{section_proof_slope_reg}}

We first present in Section \ref{section_more_general} some
\textquotedblleft structural\textquotedblright\ properties of models,
denoted \textbf{(GSA)}, that are sufficient for our needs and that are
satisfied for models of piecewise polynomial functions considered in (%
\textbf{SAPP}). Then in Sections \ref{section_proof_slope}\ and \ref%
{section_proof_hold_out}\ respectively, we prove the results stated in
Sections \ref{section_main_theorems_reg}\ and \ref{section_hold_out}, for 
\textbf{(GSA) }instead of (\textbf{SAPP}).

\subsection{A more general setting\label{section_more_general}}

\noindent \textbf{General set of assumptions: (GSA)}

\bigskip

Assume (\textbf{P1}), (\textbf{P2}), (\textbf{P3}), (\textbf{An}) and (%
\textbf{Ap}$_{u}$) of (\textbf{SAPP}). Furthermore suppose that,

\begin{description}
\item[(\textbf{Ab'})] A positive constant $A$ exists, such that for all $%
M\in \mathcal{M}_{n}$, $\left\vert Y_{i}\right\vert \leq A<\infty ,$ $%
\left\Vert s_{M}\right\Vert _{\infty }\leq A<\infty .$

\item[(\textbf{Alb})] there exists a constant $r_{\mathcal{M}}$ such that
for each $M\in \mathcal{M}_{n}$ one can find an orthonormal basis $\left(
\varphi _{k}\right) _{k=1}^{D_{M}}$ satisfying, for all $\left( \beta
_{k}\right) _{k=1}^{D_{M}}\in \mathbb{R}^{D_{M}},$%
\begin{equation*}
\left\Vert \sum_{k=1}^{D_{M}}\beta _{k}\varphi _{k}\right\Vert _{\infty
}\leq r_{\mathcal{M}}\sqrt{D_{M}}\left\vert \beta \right\vert _{\infty }%
\text{ },
\end{equation*}%
where $\left\vert \beta \right\vert _{\infty }=\max \left\{ \left\vert \beta
_{k}\right\vert ;k\in \left\{ 1,...,D_{M}\right\} \right\} $.

\item[(Ac$_{\infty }$)] a positive integer $n_{1}$ exists such that, for all 
$n\geq n_{1}$, there exist a positive constant $A_{cons}$ and an event $%
\Omega _{\infty }$ of probability at least $1-n^{-2-\alpha _{\mathcal{M}}}$,
on which for all $M\in \mathcal{M}_{n}$,%
\begin{equation}
\left\Vert s_{n}\left( M\right) -s_{M}\right\Vert _{\infty }\leq A_{cons}%
\sqrt{\frac{D_{M}\ln n}{n}}\text{ }.  \label{def_R_n_D_selection}
\end{equation}
\end{description}

\bigskip

\noindent Notice that the covariate space $\mathcal{X}$ is general in (%
\textbf{GSA}). Let us explain how assumptions (\textbf{Ab'}), (\textbf{Ad}$_{%
\leb%
}$), (\textbf{Aud}) and (\textbf{Alr}) of (\textbf{SAPP}) allow to recover (%
\textbf{Ab}), (\textbf{Alb}) and (\textbf{Ac}$_{\infty }$) of (\textbf{GSA})
in the special case of models of piecewise polynomial functions.

Assumption (\textbf{Ab'}) only differs from (\textbf{Ab}) by the fact that
the projections of the target onto the models are uniformly bounded in
sup-norm. In the general case, this is indeed not guaranteed, but
considering piecewise polynomial functions uniformly bounded in their
degrees, this follows from simple computations (see Section 5.3 in \cite%
{saum:12}). Then, assumption (\textbf{Alb}) requires the existence of a
localized orthonormal basis for each model. In the case of piecewise
polynomial functions, this is ensured by (\textbf{Ad}$_{%
\leb%
}$), (\textbf{Aud})and (\textbf{Alr}), see Lemma 8\ of \cite{saum:12}.
Finally, assumption (\textbf{Ac}$_{\infty }$) states the consistency of each
estimator for the sup-norm. Again, this is satisfied for models of piecewise
polynomial functions under assumptions (\textbf{Ad}$_{%
\leb%
}$), (\textbf{Aud}) and (\textbf{Alr}). This result is established in Lemma
9 of \cite{saum:12}.

Let us now describe a set of assumptions, less restrictive than (\textbf{%
SAPP)}, that allows to recover (\textbf{GSA}) when considering histogram
models. Lemma 5 and 6 of \cite{saum:12} allow to recover (\textbf{GSA}) from
(\textbf{SAH}) for models of histograms.

\bigskip

\noindent \textbf{Set of assumptions for histogram models: }(\textbf{SAH})

\bigskip

\noindent Given some linear histogram model $M\in \mathcal{M}_{n}$, we
denote by $\mathcal{P}_{M}$ the associated partition of $\mathcal{X}$.

\noindent Take assumptions (\textbf{P1}), (\textbf{P2}), (\textbf{P3}), (%
\textbf{An}), (\textbf{Ab}) and (\textbf{Ap}$_{u}$) from (\textbf{SAPP}).
Assume moreover,

\begin{description}
\item[(\textbf{Alrh})] there exists a positive constant $c_{\mathcal{M}%
,P}^{h}$ such that, 
\begin{equation*}
\text{for all }M\in \mathcal{M}_{n},\text{ \ }0<c_{\mathcal{M},P}^{h}\leq
\left\vert \mathcal{P}_{M}\right\vert \inf_{I\in \mathcal{P}_{M}}P^{X}\left(
I\right) \text{ .}
\end{equation*}
\end{description}

\bigskip

\noindent Theorems \ref{theorem_min_pen_reg_pp} and \ref%
{theorem_opt_pen_reg_pp} would also be valid when replacing the set of
assumptions (\textbf{SAPP}) by (\textbf{SAH}). This would lead to the
(almost exact) recovering of the assumptions and results described in
Theorems 2 and 3 of \cite{ArlotMassart:09}, concerning the selection of
least-squares estimators among histogram models.

\subsection{Proofs related to Section \protect\ref{section_main_theorems_reg}
\label{section_proof_slope}}

The following remark will be useful.

\begin{remark}
\label{remark_application_fixed_model}Since constants in (\textbf{GSA}) are
uniform over the collection $\mathcal{M}_{n}$, we deduce from Theorem 2 of 
\cite{saum:12} applied with $\alpha =2+\alpha _{\mathcal{M}}$ and $%
A_{-}=A_{+}=A_{\mathcal{M},+}$ that if assumptions (\textbf{P2}), (\textbf{%
Ab'}), (\textbf{An}), (\textbf{Alb}) and (\textbf{Ac}$_{\infty }$) hold,
then a positive constant $A_{0}$ exists, depending on $\alpha _{\mathcal{M}%
}, $ $A_{\mathcal{M},+}$ and on the constants $A,$ $\sigma _{\min }$ and $r_{%
\mathcal{M}}$ defined in (\textbf{GSA}), such that for all $M\in \mathcal{M}%
_{n}$ satisfying%
\begin{equation*}
0<A_{\mathcal{M},+}\left( \ln n\right) ^{2}\leq D_{M}\text{ },
\end{equation*}%
by setting%
\begin{equation}
\varepsilon _{n}\left( M\right) =A_{0}\max \left\{ \left( \frac{\ln n}{D_{M}}%
\right) ^{1/4};\text{ }\left( \frac{D_{M}\ln n}{n}\right) ^{1/4}\right\} 
\text{ }  \label{def_epsilon_selection}
\end{equation}%
we have, for all $n\geq n_{0}\left( A_{\mathcal{M},+},A,A_{cons},n_{1},r_{%
\mathcal{M}},\sigma _{\min },\alpha _{\mathcal{M}}\right) $, 
\begin{equation}
\mathbb{P}\left[ \left( 1-\varepsilon _{n}\left( M\right) \right) \frac{1}{4}%
\frac{D_{M}}{n}\mathcal{K}_{1,M}^{2}\leq P\left( Ks_{n}\left( M\right)
-Ks_{M}\right) \leq \left( 1+\varepsilon _{n}\left( M\right) \right) \frac{1%
}{4}\frac{D_{M}}{n}\mathcal{K}_{1,M}^{2}\right] \geq 1-10n^{-2-\alpha _{%
\mathcal{M}}}  \label{p1_selection}
\end{equation}%
and%
\begin{equation}
\mathbb{P}\left[ \left( 1-\varepsilon _{n}^{2}\left( M\right) \right) \frac{1%
}{4}\frac{D_{M}}{n}\mathcal{K}_{1,M}^{2}\leq P_{n}\left( Ks_{M}-Ks_{n}\left(
M\right) \right) \leq \left( 1+\varepsilon _{n}^{2}\left( M\right) \right) 
\frac{1}{4}\frac{D_{M}}{n}\mathcal{K}_{1,M}^{2}\right] \geq 1-5n^{-2-\alpha
_{\mathcal{M}}}  \label{p2_selection}
\end{equation}%
where $\mathcal{K}_{1,M}^{2}=1/D_{M}\sum_{k=1}^{D_{M}}\mathbb{E}\left(
\left( \psi _{1,M}\left( X,Y\right) \cdot \varphi _{k}\left( X\right)
\right) ^{2}\right) $, $\psi _{1,M}\left( X,Y\right) =-2\left( Y-s_{M}\left(
X\right) \right) ~$and $\left( \varphi _{k}\right) _{k=1}^{D_{M}}$ is an
orthonormal basis in $\left( M,\left\Vert \cdot \right\Vert _{2}\right) $.
Moreover, for all $M\in \mathcal{M}_{n}$, we have by Theorem 3 of \cite%
{saum:12}, for a positive constant $A_{u}$ depending on $A,A_{cons},r_{%
\mathcal{M}}$ and $\alpha _{\mathcal{M}}$ and for all $n\geq n_{0}\left(
A_{cons},n_{1}\right) $, 
\begin{equation}
\mathbb{P}\left[ P\left( Ks_{n}\left( M\right) -Ks_{M}\right) \geq A_{u}%
\frac{D_{M}\vee \ln n}{n}\right] \leq 3n^{-2-\alpha _{\mathcal{M}}}
\label{upper_true_selection}
\end{equation}%
and%
\begin{equation}
\mathbb{P}\left[ P_{n}\left( Ks_{M}-Ks_{n}\left( M\right) \right) \geq A_{u}%
\frac{D_{M}\vee \ln n}{n}\right] \leq 3n^{-2-\alpha _{\mathcal{M}}}\text{ }.
\label{upper_emp_selection}
\end{equation}
\end{remark}

\bigskip

Two technical lemmas are needed. In the first lemma, we intend to evaluate
the minimal penalty

\noindent $\mathbb{E}\left[ P_{n}\left( Ks_{M}-Ks_{n}\left( M\right) \right) %
\right] $ for models of dimension not too small.

\begin{lemma}
\label{mean_emp_risk}Assume (\textbf{P2}), (\textbf{Ab'}), (\textbf{An}), (%
\textbf{Alb}) and (\textbf{Ac}$_{\infty }$) of (\textbf{GSA}). Then, for
every model $M\in \mathcal{M}_{n}$ of dimension $D_{M}$ such that 
\begin{equation*}
0<A_{\mathcal{M},+}\left( \ln n\right) ^{2}\leq D_{M}\text{ ,}
\end{equation*}%
we have for all $n\geq n_{0}\left( A_{\mathcal{M},+},A,A_{cons},n_{1},r_{%
\mathcal{M}},\sigma _{\min },\alpha _{\mathcal{M}}\right) $,%
\begin{gather}
\left( 1-L_{A_{\mathcal{M},+},A,\sigma _{\min },r_{\mathcal{M}},\alpha _{%
\mathcal{M}}}\varepsilon _{n}^{2}\left( M\right) \right) \frac{D_{M}}{4n}%
\mathcal{K}_{1,M}^{2}\leq \mathbb{E}\left[ P_{n}\left( Ks_{M}-Ks_{n}\left(
M\right) \right) \right]  \label{control_mean_reg} \\
\leq \left( 1+L_{A_{\mathcal{M},+},A,\sigma _{\min },r_{\mathcal{M}},\alpha
_{\mathcal{M}}}\varepsilon _{n}^{2}\left( M\right) \right) \frac{D_{M}}{4n}%
\mathcal{K}_{1,M}^{2}\text{ ,}  \label{control_mean_2_reg}
\end{gather}%
where $\varepsilon _{n}\left( M\right) =A_{0}\max \left\{ \left( \frac{\ln n%
}{D_{M}}\right) ^{1/4};\left( \frac{D_{M}\ln n}{n}\right) ^{1/4}\right\} $
is defined in Remark \ref{remark_application_fixed_model}.
\end{lemma}

\subsection*{}%
\textbf{Proof. }As explained in Remark \ref{remark_application_fixed_model},
for all $n\geq n_{0}\left( A_{\mathcal{M},+},A,A_{cons},n_{1},r_{\mathcal{M}%
},\sigma _{\min },\alpha _{\mathcal{M}}\right) $, we thus have on an event $%
\Omega _{1}\left( M\right) $\ of probability at least $1-5n^{-2-\alpha _{%
\mathcal{M}}}$, 
\begin{equation}
\left( 1-\varepsilon _{n}\left( M\right) \right) \frac{1}{4}\frac{D_{M}}{n}%
\mathcal{K}_{1,M}^{2}\leq P_{n}\left( Ks_{M}-Ks_{n}\left( M\right) \right)
\leq \left( 1+\varepsilon _{n}\left( M\right) \right) \frac{1}{4}\frac{D_{M}%
}{n}\mathcal{K}_{1,M}^{2}\text{ ,}  \label{p2_both}
\end{equation}%
where $\varepsilon _{n}\left( M\right) =A_{0}\max \left\{ \left( \frac{\ln n%
}{D_{M}}\right) ^{1/4};\left( \frac{D_{M}\ln n}{n}\right) ^{1/4}\right\} .$
Moreover, as $\left\vert Y_{i}\right\vert \leq A$ $a.s.$ and $\left\Vert
s_{M}\right\Vert _{\infty }\leq A$ by (\textbf{Ab'}), it holds%
\begin{equation}
0\leq P_{n}\left( Ks_{M}-Ks_{n}\left( M\right) \right) \leq P_{n}Ks_{M}=%
\frac{1}{n}\sum_{i=1}^{n}\left( Y_{i}-s_{M}\left( X_{I}\right) \right)
^{2}\leq 4A^{2}  \label{majo_p2}
\end{equation}%
and as $D_{M}\geq 1$, we have%
\begin{equation}
\varepsilon _{n}\left( M\right) =A_{0}\max \left\{ \left( \frac{\ln n}{D_{M}}%
\right) ^{1/4};\left( \frac{D_{M}\ln n}{n}\right) ^{1/4}\right\} \geq
A_{0}n^{-1/8}\text{ .}  \label{mino_epsilon}
\end{equation}%
We also have%
\begin{equation}
\mathbb{E}\left[ P_{n}\left( Ks_{M}-Ks_{n}\left( M\right) \right) \right] =%
\mathbb{E}\left[ P_{n}\left( Ks_{M}-Ks_{n}\left( M\right) \right) \mathbf{1}%
_{\Omega _{1}\left( M\right) }\right] +\mathbb{E}\left[ P_{n}\left(
Ks_{M}-Ks_{n}\left( M\right) \right) \mathbf{1}_{\left( \Omega _{1}\left(
M\right) \right) ^{c}}\right] \text{ }.  \label{decompo_p2}
\end{equation}%
Now notice that by (\textbf{An}) we have $\mathcal{K}_{1,M}\geq 2\sigma
_{\min }>0$. Hence, as $D_{M}\geq 1$, it comes from (\ref{majo_p2}) and (\ref%
{mino_epsilon}) that 
\begin{equation}
0\leq \mathbb{E}\left[ P_{n}\left( Ks_{M}-Ks_{n}\left( M\right) \right) 
\mathbf{1}_{\left( \Omega _{1}\left( M\right) \right) ^{c}}\right] \leq
20A^{2}n^{-2-\alpha _{\mathcal{M}}}\leq \frac{80A^{2}}{A_{0}^{2}\sigma
_{\min }^{2}}\varepsilon _{n}^{2}\left( M\right) \frac{D_{M}}{4n}\mathcal{K}%
_{1,M}^{2}\text{ .}  \label{majo_p2_1}
\end{equation}%
Moreover, we have $\varepsilon _{n}\left( M\right) <1$ for all $n\geq
n_{0}\left( A_{0},A_{\mathcal{M},+},A_{cons}\right) $, so by (\ref{p2_both}),%
\begin{gather}
0<\left( 1-5n^{-2-\alpha _{\mathcal{M}}}\right) \left( 1-\varepsilon
_{n}^{2}\left( M\right) \right) \frac{D_{M}}{4n}\mathcal{K}_{1,M}^{2}\leq 
\mathbb{E}\left[ P_{n}\left( Ks_{M}-Ks_{n}\left( M\right) \right) \mathbf{1}%
_{\Omega _{1}\left( M\right) }\right]  \label{mino_p2_2} \\
\leq \left( 1+\varepsilon _{n}^{2}\left( M\right) \right) \frac{D_{M}}{4n}%
\mathcal{K}_{1,M}^{2}\text{ }.  \label{majo_p2_2}
\end{gather}%
Finally, noticing that $n^{-2-\alpha _{\mathcal{M}}}\leq
A_{0}^{-2}\varepsilon _{n}^{2}\left( M\right) $ by (\ref{mino_epsilon}), we
use (\ref{majo_p2_1}), (\ref{mino_p2_2}) and (\ref{majo_p2_2}) in (\ref%
{decompo_p2}) to conclude by straightforward computations that 
\begin{equation*}
L_{A_{\mathcal{M},+},A,\sigma _{\min },r_{\mathcal{M}},\alpha _{\mathcal{M}%
}}=\frac{80A^{2}}{A_{0}^{2}\sigma _{\min }^{2}}+5A_{0}^{-2}+1
\end{equation*}%
is convenient in (\ref{control_mean_reg}) and (\ref{control_mean_2_reg}), as 
$A_{0}$ only depends on $\alpha _{\mathcal{M}},$ $A_{\mathcal{M},+},$ $A,$ $%
\sigma _{\min }$ and $r_{\mathcal{M}}$.\ 
$\blacksquare$%

\begin{lemma}
\label{delta_bar}Let $\alpha >0$. Assume that (\textbf{Ab'}) of (\textbf{GSA}%
) is satisfied. Then there exists a positive constant $A_{d}$, depending
only in $A,$ $A_{\mathcal{M},+},$ $\sigma _{\min }$ and $\alpha $ such that,
by setting $\bar{\delta}\left( M\right) =\left( P_{n}-P\right) \left(
Ks_{M}-Ks_{\ast }\right) $, we have for all $M\in \mathcal{M}_{n}$, 
\begin{equation}
\mathbb{P}\left( \left\vert \bar{\delta}\left( M\right) \right\vert \geq
A_{d}\left( \sqrt{\frac{\ell \left( s_{\ast },s_{M}\right) \ln n}{n}}+\frac{%
\ln n}{n}\right) \right) \leq 2n^{-\alpha }\text{ }.  \label{delta_bar_small}
\end{equation}%
If moreover, assumptions (\textbf{P2}), (\textbf{An}), (\textbf{Alb}) and (%
\textbf{Ac}$_{\infty }$) of (\textbf{GSA}) hold, then for all $M\in \mathcal{%
M}_{n}$ such that $A_{\mathcal{M},+}\left( \ln n\right) ^{2}\leq D_{M}$ and
for all $n\geq n_{0}\left( A_{\mathcal{M},+},A,A_{cons},n_{1},r_{\mathcal{M}%
},\sigma _{\min },\alpha \right) $, we have%
\begin{equation}
\mathbb{P}\left( \left\vert \bar{\delta}\left( M\right) \right\vert \geq 
\frac{\ell \left( s_{\ast },s_{M}\right) }{\sqrt{D_{M}}}+A_{d}\frac{\ln n}{%
\sqrt{D_{M}}}\mathbb{E}\left[ 
\p_{2}\left(M\right)%
\right] \right) \leq 2n^{-\alpha }\text{ ,}  \label{delta_bar_reasonable}
\end{equation}%
where $%
\p_{2}\left(M\right)%
:=P_{n}\left( Ks_{M}-Ks_{n}\left( M\right) \right) \geq 0$.
\end{lemma}

\subsection*{}%
\textbf{Proof. }We set%
\begin{equation}
A_{d}=\max \left\{ 4A\sqrt{\alpha };\text{ }\frac{8A^{2}}{3}\alpha ;\text{ }%
\frac{8A^{2}\alpha }{\sqrt{A_{\mathcal{M},+}}\sigma _{\min }^{2}}+\frac{%
16A^{2}\alpha }{3A_{\mathcal{M},+}\sigma _{\min }}\right\} \text{ .}
\label{calcul_Ad}
\end{equation}%
Since by (\textbf{Ab'}) we have $\left\vert Y\right\vert \leq A$ $a.s.$ and $%
\left\Vert s_{M}\right\Vert _{\infty }\leq A$, it holds $\left\Vert s_{\ast
}\right\Vert _{\infty }=\left\Vert \mathbb{E}\left[ Y\left\vert X\right. %
\right] \right\Vert _{\infty }\leq A$, and so $\left\Vert s_{M}-s_{\ast
}\right\Vert _{\infty }\leq 2A.$ Next, we apply Bernstein's inequality (see
Proposition 2.9 of \cite{Massart:07}) to $\bar{\delta}\left( M\right)
=\left( P_{n}-P\right) \left( Ks_{M}-Ks_{\ast }\right) .$ Notice that%
\begin{equation*}
K\left( s_{M}\right) \left( x,y\right) -K\left( s_{\ast }\right) \left(
x,y\right) =\left( s_{M}\left( x\right) -s_{\ast }\left( x\right) \right)
\left( s_{M}\left( x\right) +s_{\ast }\left( x\right) -2y\right) \text{ },
\end{equation*}%
hence $\left\Vert Ks_{M}-Ks_{\ast }\right\Vert _{\infty }\leq 8A^{2}.$
Moreover, as $\mathbb{E}\left[ Y-s_{\ast }\left( X\right) \left\vert
X\right. \right] =0$ and $\mathbb{E}\left[ \left( Y-s_{\ast }\left( X\right)
\right) ^{2}\left\vert X\right. \right] \leq \frac{\left( 2A\right) ^{2}}{4}%
=A^{2}$ we have%
\begin{align*}
& \mathbb{E}\left[ \left( Ks_{M}\left( X,Y\right) -Ks_{\ast }\left(
X,Y\right) \right) ^{2}\right] \\
& =\mathbb{E}\left[ \left( 4\left( Y-s_{\ast }\left( X\right) \right)
^{2}+\left( s_{M}\left( X\right) -s_{\ast }\left( X\right) \right)
^{2}\right) \left( s_{M}\left( X\right) -s_{\ast }\left( X\right) \right)
^{2}\right] \\
& \leq 8A^{2}\mathbb{E}\left[ \left( s_{M}\left( X\right) -s_{\ast }\left(
X\right) \right) ^{2}\right] =8A^{2}\ell \left( s_{\ast },s_{M}\right) ,
\end{align*}%
and therefore, by Bernstein's inequality we have for all $x>0,$%
\begin{equation*}
\mathbb{P}\left( \left\vert \bar{\delta}\left( M\right) \right\vert \geq 
\sqrt{\frac{16A^{2}\ell \left( s_{\ast },s_{M}\right) x}{n}}+\frac{8A^{2}x}{%
3n}\right) \leq 2\exp \left( -x\right) \text{ }.
\end{equation*}%
By taking $x=\alpha \ln n$, we then have%
\begin{equation}
\mathbb{P}\left( \left\vert \bar{\delta}\left( M\right) \right\vert \geq 
\sqrt{\frac{16A^{2}\alpha \ell \left( s_{\ast },s_{M}\right) \ln n}{n}}+%
\frac{8A^{2}\alpha \ln n}{3n}\right) \leq 2n^{-\alpha }\text{ ,}
\label{delta_small_models}
\end{equation}%
which gives the first part of Lemma \ref{delta_bar} for $A_{d}$ given in (%
\ref{calcul_Ad}). Now, by noticing the fact that $2\sqrt{ab}\leq a\eta
+b\eta ^{-1}$ for all $\eta >0$, and using it in (\ref{delta_small_models})
with $a=\ell \left( s_{\ast },s_{M}\right) $, $b=\frac{4A^{2}\alpha \ln n}{n}
$ and $\eta =D_{M}^{-1/2}$ , we obtain%
\begin{equation}
\mathbb{P}\left( \left\vert \bar{\delta}\left( M\right) \right\vert \geq 
\frac{\ell \left( s_{\ast },s_{M}\right) }{\sqrt{D_{M}}}+\left( 4\sqrt{D_{M}}%
+\frac{8}{3}\right) \frac{A^{2}\alpha \ln n}{n}\right) \leq 2n^{-\alpha }%
\text{ }.  \label{delta_2}
\end{equation}%
Then, for a model $M\in \mathcal{M}_{n}$ such that $A_{\mathcal{M},+}\left(
\ln n\right) ^{2}\leq D_{M}$, we apply Lemma \ref{mean_emp_risk} and by (\ref%
{control_mean_reg}), it holds for all $n\geq n_{0}\left( A_{\mathcal{M}%
,+},A,A_{cons},n_{1},r_{\mathcal{M}},\sigma _{\min },\alpha _{\mathcal{M}%
}\right) $,%
\begin{equation}
\left( 1-L_{A_{\mathcal{M},-},A,\sigma _{\min },r_{\mathcal{M}},\alpha _{%
\mathcal{M}}}\varepsilon _{n}^{2}\left( M\right) \right) \frac{D_{M}}{4n}%
\mathcal{K}_{1,M}^{2}\leq \mathbb{E}\left[ 
\p_{2}\left(M\right)%
\right]  \label{mino_mean_2}
\end{equation}%
where $\varepsilon _{n}\left( M\right) =A_{0}\max \left\{ \left( \frac{\ln n%
}{D_{M}}\right) ^{1/4};\left( \frac{D_{M}\ln n}{n}\right) ^{1/4}\right\} $.
Moreover, as $D_{M}\leq A_{\mathcal{M},+}n\left( \ln n\right) ^{-2}$ by (%
\textbf{P2}) and $A_{\mathcal{M},+}\left( \ln n\right) ^{2}\leq D_{M}$, we
deduce that for all $n\geq n_{0}\left( A_{\mathcal{M},+},A,A_{cons},r_{%
\mathcal{M}},\sigma _{\min },\alpha _{\mathcal{M}}\right) $, 
\begin{equation*}
L_{A_{\mathcal{M},-},A,\sigma _{\min },r_{\mathcal{M}},\alpha _{\mathcal{M}%
}}\varepsilon _{n}^{2}\left( M\right) \leq 1/2\text{ .}
\end{equation*}%
Now, since $\mathcal{K}_{1,M}\geq 2\sigma _{\min }>0$ by (\textbf{An}), we
have by (\ref{mino_mean_2}), $\mathbb{E}\left[ 
\p_{2}\left(M\right)%
\right] \geq \frac{\sigma _{\min }^{2}}{2}\frac{D_{M}}{n}$ for all

\noindent $n\geq n_{0}\left( A_{\mathcal{M},+},A,A_{cons},n_{1},r_{\mathcal{M%
}},\sigma _{\min },\alpha _{\mathcal{M}}\right) $. This allows, using (\ref%
{delta_2}), to conclude the proof for the value of $A_{d}$ given in (\ref%
{calcul_Ad}) by simple computations.\ 
$\blacksquare$%

\bigskip

\noindent In order to avoid cumbersome notations in the proofs of Theorems %
\ref{theorem_opt_pen_reg_pp} and \ref{theorem_min_pen_reg_pp}, when generic
constants $L$ and $n_{0}$ depend on constants defined in the general set of
assumptions stated in Section \ref{section_more_general}, we will note $L_{%
\text{(\textbf{GSA})}}$ and $n_{0}\left( \text{\textbf{GSA}}\right) $. The
values of these constants may change from line to line.

\bigskip

\noindent \textbf{Proof of Theorem \ref{theorem_opt_pen_reg_pp}. }From the
definition of the selected model $%
\widehat{M}%
$ given in (\ref{def_proc_2_reg}), $%
\widehat{M}%
$ minimizes 
\begin{equation}
\crit%
\left( M\right) :=P_{n}\left( Ks_{n}\left( M\right) \right) +%
\pen%
\left( M\right) \text{ ,}  \label{def_crit}
\end{equation}%
over the models $M\in \mathcal{M}_{n}$. Hence, $%
\widehat{M}%
$ also minimizes%
\begin{equation}
\crit%
^{\prime }\left( M\right) :=%
\crit%
\left( M\right) -P_{n}\left( Ks_{\ast }\right) \text{ ,}
\label{def_criterion_prime_1}
\end{equation}%
over the collection $\mathcal{M}_{n}$. Let us write 
\begin{align*}
\ell \left( s_{\ast },s_{n}\left( M\right) \right) & =P\left( Ks_{n}\left(
M\right) -Ks_{\ast }\right) \\
& =P_{n}\left( Ks_{n}\left( M\right) \right) +P_{n}\left(
Ks_{M}-Ks_{n}\left( M\right) \right) +\left( P_{n}-P\right) \left( Ks_{\ast
}-Ks_{M}\right) \\
& +P\left( Ks_{n}\left( M\right) -Ks_{M}\right) -P_{n}\left( Ks_{\ast
}\right) \text{ }.
\end{align*}%
By setting 
\begin{equation*}
\p_{1}\left(M\right)%
=P\left( Ks_{n}\left( M\right) -Ks_{M}\right) \text{ ,}
\end{equation*}%
\begin{equation*}
\p_{2}\left(M\right)%
=P_{n}\left( Ks_{M}-Ks_{n}\left( M\right) \right) \text{ ,}
\end{equation*}%
\begin{equation*}
\bar{\delta}\left( M\right) =\left( P_{n}-P\right) \left( Ks_{M}-Ks_{\ast
}\right)
\end{equation*}%
and%
\begin{equation*}
\pen%
_{\text{id}}^{\prime }\left( M\right) =%
\p_{1}\left(M\right)%
+%
\p_{2}\left(M\right)%
-\bar{\delta}\left( M\right) \text{ ,}
\end{equation*}%
we have%
\begin{equation}
\ell \left( s_{\ast },s_{n}\left( M\right) \right) =P_{n}\left( Ks_{n}\left(
M\right) \right) +%
\p_{1}\left(M\right)%
+%
\p_{2}\left(M\right)%
-\bar{\delta}\left( M\right) -P_{n}\left( Ks_{\ast }\right)
\label{decompo_excess}
\end{equation}%
and by (\ref{def_criterion_prime_1}),%
\begin{equation}
\crit%
^{\prime }\left( M\right) =\ell \left( s_{\ast },s_{n}\left( M\right)
\right) +\left( 
\pen%
\left( M\right) -%
\pen%
_{\text{id}}^{\prime }\left( M\right) \right) \text{ .}
\label{formule_crit_prime}
\end{equation}%
As $%
\widehat{M}%
$ minimizes $%
\crit%
^{\prime }$ over $\mathcal{M}_{n}$, it is therefore sufficient by (\ref%
{formule_crit_prime}), to control $%
\pen%
\left( M\right) -%
\pen%
_{\text{id}}^{\prime }\left( M\right) $ - or equivalently $%
\crit%
^{\prime }\left( M\right) $ - in terms of the excess loss $\ell \left(
s_{\ast },s_{n}\left( M\right) \right) $, for every $M\in \mathcal{M}_{n}$,
in order to derive oracle inequalities. Let $\Omega _{n}$ be the event on
which:

\begin{itemize}
\item For all models $M\in \mathcal{M}_{n}$ of dimension $D_{M}$ such that $%
A_{\mathcal{M},+}\left( \ln n\right) ^{3}\leq D_{M}$, (\ref{pen_id_pp})
holds and 
\begin{align}
\left\vert 
\p_{1}\left(M\right)%
-\mathbb{E}\left[ 
\p_{2}\left(M\right)%
\right] \right\vert & \leq L_{\text{(\textbf{GSA})}}\varepsilon _{n}\left(
M\right) \mathbb{E}\left[ 
\p_{2}\left(M\right)%
\right]  \label{line_1} \\
\left\vert 
\p_{2}\left(M\right)%
-\mathbb{E}\left[ 
\p_{2}\left(M\right)%
\right] \right\vert & \leq L_{\text{(\textbf{GSA})}}\varepsilon
_{n}^{2}\left( M\right) \mathbb{E}\left[ 
\p_{2}\left(M\right)%
\right]  \label{line_2} \\
\left\vert \bar{\delta}\left( M\right) \right\vert & \leq \frac{\ell \left(
s_{\ast },s_{M}\right) }{\sqrt{D_{M}}}+L_{\text{(\textbf{GSA})}}\frac{\ln n}{%
\sqrt{D_{M}}}\mathbb{E}\left[ 
\p_{2}\left(M\right)%
\right]  \label{line_3_ter} \\
\left\vert \bar{\delta}\left( M\right) \right\vert & \leq L_{\text{(\textbf{%
GSA})}}\left( \sqrt{\frac{\ell \left( s_{\ast },s_{M}\right) \ln n}{n}}+%
\frac{\ln n}{n}\right)  \label{line_4}
\end{align}

\item For all models $M\in \mathcal{M}_{n}$ of dimension $D_{M}$ such that $%
D_{M}\leq A_{\mathcal{M},+}\left( \ln n\right) ^{3}$, (\ref{pen_id_2_pp})
holds together with%
\begin{align}
\left\vert \bar{\delta}\left( M\right) \right\vert & \leq L_{\text{(\textbf{%
GSA})}}\left( \sqrt{\frac{\ell \left( s_{\ast },s_{M}\right) \ln n}{n}}+%
\frac{\ln n}{n}\right)  \label{delta_bar_small_models} \\
\p_{2}\left(M\right)%
& \leq L_{\text{(\textbf{GSA})}}\frac{D_{M}\vee \ln n}{n}\leq L_{\text{(%
\textbf{GSA})}}\frac{\left( \ln n\right) ^{3}}{n}  \label{p_2_small_models_1}
\\
\p_{1}\left(M\right)%
& \leq L_{\text{(\textbf{GSA})}}\frac{D_{M}\vee \ln n}{n}\leq L_{\text{(%
\textbf{GSA})}}\frac{\left( \ln n\right) ^{3}}{n}  \label{p_1_small_models}
\end{align}
\end{itemize}

\noindent By (\ref{p1_selection}), (\ref{p2_selection}), (\ref%
{upper_true_selection}) and (\ref{upper_emp_selection}) in Remark \ref%
{remark_application_fixed_model}, Lemma \ref{mean_emp_risk}, Lemma \ref%
{delta_bar} applied with $\alpha =2+\alpha _{\mathcal{M}}$, and since (\ref%
{pen_id_pp}) holds with probability at least $1-A_{p}n^{-2}$, we get for all 
$n\geq n_{0}\left( \text{\textbf{GSA}}\right) $, 
\begin{equation*}
\mathbb{P}\left( \Omega _{n}\right) \geq 1-A_{p}n^{-2}-24\sum_{M\in \mathcal{%
M}_{n}}n^{-2-\alpha _{\mathcal{M}}}\geq 1-L_{A_{p},c_{\mathcal{M}}}n^{-2}%
\text{ }.
\end{equation*}

\bigskip

\noindent \underline{\textbf{Control on the criterion }$%
\crit%
^{\prime }$\textbf{\ for models of dimension not too small: }}

\bigskip

\noindent We consider models $M\in \mathcal{M}_{n}$ such that $A_{\mathcal{M}%
,+}\left( \ln n\right) ^{3}\leq D_{M}$. Notice that (\ref{line_3_ter})
implies by (\ref{def_epsilon_selection}) that, for all $M\in \mathcal{M}_{n}$
such that $A_{\mathcal{M},+}\left( \ln n\right) ^{3}\leq D_{M}$, for all $%
n\geq n_{0}\left( \text{\textbf{GSA}}\right) $, 
\begin{eqnarray*}
\left\vert \bar{\delta}\left( M\right) \right\vert &\leq &L_{\text{(\textbf{%
GSA})}}\left( \frac{\left( \ln n\right) ^{3}}{D_{M}}%
\cdot%
\frac{\ln n}{D_{M}}\right) ^{1/4}\times \mathbb{E}\left[ \ell \left( s_{\ast
},s_{M}\right) +%
\p_{2}\left(M\right)%
\right] \\
&\leq &L_{\text{(\textbf{GSA})}}\varepsilon _{n}\left( M\right) \mathbb{E}%
\left[ \ell \left( s_{\ast },s_{M}\right) +%
\p_{2}\left(M\right)%
\right] \text{ ,}
\end{eqnarray*}%
so that on $\Omega _{n}$ we have, for all models $M\in \mathcal{M}_{n}$ such
that $A_{\mathcal{M},+}\left( \ln n\right) ^{3}\leq D_{M}$, 
\begin{align}
& \left\vert 
\pen%
_{\text{id}}^{\prime }\left( M\right) -%
\pen%
\left( M\right) \right\vert  \notag \\
& \leq \left\vert 
\p_{1}\left(M\right)%
+%
\p_{2}\left(M\right)%
-%
\pen%
\left( M\right) \right\vert +\left\vert \bar{\delta}\left( M\right)
\right\vert  \notag \\
& \leq \left\vert 
\p_{1}\left(M\right)%
+%
\p_{2}\left(M\right)%
-2\mathbb{E}\left[ 
\p_{2}\left(M\right)%
\right] \right\vert +\left( L_{\text{(\textbf{GSA})}}\varepsilon _{n}\left(
M\right) +\delta \right) \mathbb{E}\left[ \ell \left( s_{\ast },s_{M}\right)
+%
\p_{2}\left(M\right)%
\right]  \notag \\
& \leq \left( \delta +L_{\text{(\textbf{GSA})}}\varepsilon _{n}\left(
M\right) \right) \mathbb{E}\left[ \ell \left( s_{\ast },s_{M}\right) +%
\p_{2}\left(M\right)%
\right] \text{ }.  \label{pen-pen_id}
\end{align}%
Now notice that using (\textbf{P2}) in (\ref{def_epsilon_selection}) gives
that for all models $M\in \mathcal{M}_{n}$ such that $A_{\mathcal{M}%
,+}\left( \ln n\right) ^{3}\leq D_{M}$ and for all $n\geq n_{0}\left( \text{%
\textbf{GSA}}\right) $, $0<L_{\text{(\textbf{GSA})}}\varepsilon _{n}\left(
M\right) \leq \frac{1}{2}$. As $\ell \left( s_{\ast },s_{n}\left( M\right)
\right) =\ell \left( s_{\ast },s_{M}\right) +%
\p_{1}\left(M\right)%
$, we thus have on $\Omega _{n}$, for all $n\geq n_{0}\left( \text{\textbf{%
GSA}}\right) $,%
\begin{align}
0& \leq \mathbb{E}\left[ \ell \left( s_{\ast },s_{M}\right) +%
\p_{2}\left(M\right)%
\right]  \notag \\
& \leq \ell \left( s_{\ast },s_{n}\left( M\right) \right) +\left\vert 
\p_{1}\left(M\right)%
-\mathbb{E}\left[ 
\p_{2}\left(M\right)%
\right] \right\vert  \notag \\
& \leq \ell \left( s_{\ast },s_{n}\left( M\right) \right) +\frac{L_{\text{(%
\textbf{GSA})}}\varepsilon _{n}\left( M\right) }{1-L_{\text{(\textbf{GSA})}%
}\varepsilon _{n}\left( M\right) }%
\p_{1}\left(M\right)%
\text{ \ \ \ by (\ref{line_1})}  \notag \\
& \leq \frac{1+L_{\text{(\textbf{GSA})}}\varepsilon _{n}\left( M\right) }{%
1-L_{\text{(\textbf{GSA})}}\varepsilon _{n}\left( M\right) }\ell \left(
s_{\ast },s_{n}\left( M\right) \right)  \notag \\
& \leq \left( 1+L_{\text{(\textbf{GSA})}}\varepsilon _{n}\left( M\right)
\right) \ell \left( s_{\ast },s_{n}\left( M\right) \right) \text{ .}
\label{control_l}
\end{align}%
Hence, using (\ref{control_l}) in (\ref{pen-pen_id}), we have on $\Omega
_{n} $ for all models $M\in \mathcal{M}_{n}$ such that $A_{\mathcal{M}%
,+}\left( \ln n\right) ^{3}\leq D_{M}$ and for all $n\geq n_{0}\left( \text{%
\textbf{GSA}}\right) $,%
\begin{equation}
\left\vert 
\pen%
_{\text{id}}^{\prime }\left( M\right) -%
\pen%
\left( M\right) \right\vert \leq \left( \delta +L_{\text{(\textbf{GSA})}%
}\varepsilon _{n}\left( M\right) \right) \ell \left( s_{\ast },s_{n}\left(
M\right) \right) \text{ .}  \label{pen-pen_id_2}
\end{equation}%
Consequently, for all models $M\in \mathcal{M}_{n}$ such that $A_{\mathcal{M}%
,+}\left( \ln n\right) ^{3}\leq D_{M}$ and for all $n\geq n_{0}\left( \text{%
\textbf{GSA}}\right) $, it holds on $\Omega _{n}$, using (\ref%
{formule_crit_prime}) and (\ref{pen-pen_id_2}),%
\begin{equation}
\left( 1-\delta -L_{\text{(\textbf{GSA})}}\varepsilon _{n}\left( M\right)
\right) \ell \left( s_{\ast },s_{n}\left( M\right) \right) \leq 
\crit%
^{\prime }\left( M\right) \leq \left( 1+\delta +L_{\text{(\textbf{GSA})}%
}\varepsilon _{n}\left( M\right) \right) \ell \left( s_{\ast },s_{n}\left(
M\right) \right) \text{ .}  \label{control_crit_prime_large}
\end{equation}

\bigskip

\noindent \underline{\textbf{Control on the criterion }$%
\crit%
^{\prime }$\textbf{\ for models of small dimension: }}

\bigskip

\noindent We consider models $M\in \mathcal{M}_{n}$ such that $D_{M}\leq A_{%
\mathcal{M},+}\left( \ln n\right) ^{3}$. By (\ref{pen_id_2_pp}), (\ref%
{delta_bar_small_models}) and (\ref{p_2_small_models_1}), it holds on $%
\Omega _{n}$, for any $\tau >0$ and for all $M\in \mathcal{M}_{n}$ such that 
$D_{M}\leq A_{\mathcal{M},+}\left( \ln n\right) ^{3}$,

\begin{align}
& \left\vert 
\pen%
_{\text{id}}^{\prime }\left( M\right) -%
\pen%
\left( M\right) \right\vert  \notag \\
& \leq 
\p_{1}\left(M\right)%
+%
\p_{2}\left(M\right)%
+\left\vert 
\pen%
\left( M\right) \right\vert +\left\vert \bar{\delta}\left( M\right)
\right\vert  \notag \\
& \leq L_{\text{(\textbf{GSA})}}\frac{\left( \ln n\right) ^{3}}{n}+A_{r}%
\frac{\ell \left( s_{\ast },s_{M}\right) }{\left( \ln n\right) ^{2}}+A_{r}%
\frac{\left( \ln n\right) ^{3}}{n}+L_{\text{(\textbf{GSA})}}\left( \sqrt{%
\frac{\ell \left( s_{\ast },s_{M}\right) \ln n}{n}}+\frac{\ln n}{n}\right) 
\notag \\
& \leq L_{\text{(\textbf{GSA}),}A_{r}}\left( \frac{\left( \ln n\right) ^{3}}{%
n}+\frac{\ell \left( s_{\ast },s_{M}\right) }{\left( \ln n\right) ^{2}}%
\right) +\tau \ell \left( s_{\ast },s_{M}\right) +\left( \tau ^{-1}+1\right)
L_{\text{(\textbf{GSA})}}\frac{\ln n}{n}  \notag \\
& \leq L_{\text{(\textbf{GSA}),}A_{r}}\left( \frac{\left( \ln n\right) ^{3}}{%
n}+\frac{\ell \left( s_{\ast },s_{M}\right) }{\left( \ln n\right) ^{2}}%
\right) +\tau \ell \left( s_{\ast },s_{n}\left( M\right) \right) +\left(
\tau ^{-1}+1\right) L_{\text{(\textbf{GSA})}}\frac{\ln n}{n}\text{ .}
\label{majo_small}
\end{align}%
Hence, by taking $\tau =\left( \ln n\right) ^{-2}$ in (\ref{majo_small}) we
get that for all $M\in \mathcal{M}_{n}$ such that $D_{M}\leq A_{\mathcal{M}%
,+}\left( \ln n\right) ^{3}$, it holds on $\Omega _{n}$,%
\begin{equation}
\left\vert 
\pen%
_{\text{id}}^{\prime }\left( M\right) -%
\pen%
\left( M\right) \right\vert \leq L_{\text{(\textbf{GSA}),}A_{r}}\left( \frac{%
\ell \left( s_{\ast },s_{n}\left( M\right) \right) }{\left( \ln n\right) ^{2}%
}+\frac{\left( \ln n\right) ^{3}}{n}\right) \text{ .}  \label{majo_small_2}
\end{equation}%
Moreover, by (\ref{formule_crit_prime}) and (\ref{majo_small_2}), we have on
the event $\Omega _{n}$, for all $M\in \mathcal{M}_{n}$ such that $D_{M}\leq
A_{\mathcal{M},+}\left( \ln n\right) ^{3}$,%
\begin{gather}
\left( 1-L_{\text{(\textbf{GSA}),}A_{r}}\left( \ln n\right) ^{-2}\right)
\ell \left( s_{\ast },s_{n}\left( M\right) \right) -L_{\text{(\textbf{GSA}),}%
A_{r}}\frac{\left( \ln n\right) ^{3}}{n}\leq 
\crit%
^{\prime }\left( M\right)  \label{lower_crit_prime_small} \\
\leq \left( 1+L_{\text{(\textbf{GSA}),}A_{r}}\left( \ln n\right)
^{-2}\right) \ell \left( s_{\ast },s_{n}\left( M\right) \right) +L_{\text{(%
\textbf{GSA}),}A_{r}}\frac{\left( \ln n\right) ^{3}}{n}\text{ .}
\label{upper_crit_prime_small}
\end{gather}

\bigskip

\noindent \underline{\textbf{Oracle inequalities: }}

\bigskip

\noindent Recall that by the definition given in (\ref{oracle_model}), an
oracle model satisfies 
\begin{equation}
M_{\ast }\in \arg \min_{M\in \mathcal{M}_{n}}\left\{ \ell \left( s_{\ast
},s_{n}\left( M\right) \right) \right\} \text{ }.  \label{def_Mstar_conclu}
\end{equation}%
By Lemmas \ref{lemma_selected_model}\ and \ref{lemma_oracle_model}\ below,
we control on $\Omega _{n}$ the dimensions of the selected model $%
\widehat{M}%
$ and the oracle model $M_{\ast }$. More precisely, by (\ref%
{control_dim_hat_gene}) and (\ref{control_dim_star_gene}), we have on $%
\Omega _{n}$, for any $\eta \in \left( 0,\beta _{+}/\left( 1+\beta
_{+}\right) \right) $ and for all $n\geq n_{0}\left( \text{(\textbf{GSA})}%
,\eta ,\delta \right) $, 
\begin{align}
D_{%
\widehat{M}%
}& \leq n^{1/\left( 1+\beta _{+}\right) +\eta }\text{ ,}
\label{local_dim_hat_gene} \\
D_{M_{\ast }}& \leq n^{1/\left( 1+\beta _{+}\right) +\eta }\text{ }.
\label{local_dim_star_gene}
\end{align}%
Now, from (\ref{local_dim_hat_gene}) we distinguish two cases in order to
control $%
\crit%
^{\prime }\left( 
\widehat{M}%
\right) $. If $A_{\mathcal{M},+}\left( \ln n\right) ^{3}\leq D_{%
\widehat{M}%
}\leq n^{1/\left( 1+\beta _{+}\right) +\eta }$, we get by (\ref%
{control_crit_prime_large}), for all $n\geq n_{0}\left( \text{\textbf{GSA}}%
\right) $,%
\begin{equation}
\crit%
^{\prime }\left( 
\widehat{M}%
\right) \geq \left( 1-\delta -L_{\text{(\textbf{GSA})}}\varepsilon
_{n}\left( 
\widehat{M}%
\right) \right) \ell \left( s_{\ast },s_{n}\left( 
\widehat{M}%
\right) \right) \text{ .}  \label{lower_crit_hat_large}
\end{equation}%
Otherwise, if $D_{%
\widehat{M}%
}\leq A_{\mathcal{M},+}\left( \ln n\right) ^{3}$, we get by (\ref%
{lower_crit_prime_small}), 
\begin{equation}
\left( 1-L_{\text{(\textbf{GSA}),}A_{r}}\left( \ln n\right) ^{-2}\right)
\ell \left( s_{\ast },s_{n}\left( 
\widehat{M}%
\right) \right) -L_{\text{(\textbf{GSA}),}A_{r}}\frac{\left( \ln n\right)
^{3}}{n}\leq 
\crit%
^{\prime }\left( 
\widehat{M}%
\right) \text{ .}  \label{lower_crit_hat_small}
\end{equation}%
Let us denote $\mathcal{S}_{n}=\left\{ M\in \mathcal{M}_{n};\text{ }A_{%
\mathcal{M},+}\left( \ln n\right) ^{3}\leq D_{M}\leq n^{1/\left( 1+\beta
_{+}\right) +\eta }\right\} $. In all cases, we have by (\ref%
{lower_crit_hat_large}) and (\ref{lower_crit_hat_small}), for all $n\geq
n_{0}\left( \text{\textbf{GSA}}\right) $,%
\begin{gather}
\crit%
^{\prime }\left( 
\widehat{M}%
\right) \geq \left( 1-\delta -L_{\text{(\textbf{GSA}),}A_{r}}\left( \left(
\ln n\right) ^{-2}+\sup_{M\in \mathcal{S}_{n}}\varepsilon _{n}\left(
M\right) \right) \right) \ell \left( s_{\ast },s_{n}\left( 
\widehat{M}%
\right) \right)  \notag \\
-L_{\text{(\textbf{GSA}),}A_{r}}\frac{\left( \ln n\right) ^{3}}{n}\text{ .}
\label{lower_crit_hat_gene}
\end{gather}%
Similarly, from (\ref{local_dim_star_gene}) we distinguish two cases in
order to control $%
\crit%
^{\prime }\left( M_{\ast }\right) $. If $A_{\mathcal{M},+}\left( \ln
n\right) ^{3}\leq D_{M_{\ast }}\leq n^{1/\left( 1+\beta _{+}\right) +\eta }$%
, we get by (\ref{control_crit_prime_large}), for all $n\geq n_{0}\left( 
\text{\textbf{GSA}}\right) $,%
\begin{equation}
\crit%
^{\prime }\left( M_{\ast }\right) \leq \left( 1+\delta +L_{\text{(\textbf{GSA%
})}}\varepsilon _{n}\left( M_{\ast }\right) \right) \ell \left( s_{\ast
},s_{n}\left( M_{\ast }\right) \right) \text{ .}
\label{upper_crit_star_large}
\end{equation}%
Otherwise, if $D_{M_{\ast }}\leq A_{\mathcal{M},+}\left( \ln n\right) ^{3}$,
we get by (\ref{upper_crit_prime_small}), 
\begin{equation}
\crit%
^{\prime }\left( M_{\ast }\right) \leq \left( 1+L_{\text{(\textbf{GSA}),}%
A_{r}}\left( \ln n\right) ^{-2}\right) \ell \left( s_{\ast },s_{n}\left(
M_{\ast }\right) \right) +L_{\text{(\textbf{GSA}),}A_{r}}\frac{\left( \ln
n\right) ^{3}}{n}\text{ .}  \label{upper_crit_star_small}
\end{equation}%
In all cases, we deduce from (\ref{upper_crit_star_large}) and (\ref%
{upper_crit_star_small}) that we have for all $n\geq n_{0}\left( \text{(%
\textbf{GSA}),}\delta \right) $,%
\begin{gather}
\crit%
^{\prime }\left( M_{\ast }\right) \leq \left( 1+\delta +L_{\text{(\textbf{GSA%
}),}A_{r}}\left( \left( \ln n\right) ^{-2}+\sup_{M\in \mathcal{S}%
_{n}}\varepsilon _{n}\left( M\right) \right) \right) \ell \left( s_{\ast
},s_{n}\left( M_{\ast }\right) \right)  \notag \\
+L_{\text{(\textbf{GSA}),}A_{r}}\frac{\left( \ln n\right) ^{3}}{n}\text{ .}
\label{upper_crit_star_gene}
\end{gather}%
Hence, by setting 
\begin{equation*}
\theta _{n}=L_{\text{(\textbf{GSA}),}A_{r}}\left( \left( \ln n\right)
^{-2}+\sup_{M\in \mathcal{S}_{n}}\varepsilon _{n}\left( M\right) \right) 
\text{ ,}
\end{equation*}%
we have by (\ref{def_epsilon_selection}), for all $n\geq n_{0}\left( \text{(%
\textbf{GSA})},\eta ,\delta \right) $,%
\begin{equation*}
\theta _{n}\leq \frac{L_{\text{(\textbf{GSA}),}A_{r}}}{\left( \ln n\right)
^{1/4}}\text{ , \ \ \ }\theta _{n}<\frac{1-\delta }{2}
\end{equation*}%
and we deduce from (\ref{lower_crit_hat_gene}) and (\ref%
{upper_crit_star_gene}), since $\frac{1}{1-x}\leq 1+2x$ for all $x\in \left[
0,\frac{1}{2}\right) $, that for all $n\geq n_{0}\left( \text{(\textbf{GSA})}%
,\eta ,\delta \right) $, it holds on $\Omega _{n}$,%
\begin{align}
\ell \left( s_{\ast },s_{n}\left( 
\widehat{M}%
\right) \right) & \leq \left( \frac{1+\delta +\theta _{n}}{1-\delta -\theta
_{n}}\right) \ell \left( s_{\ast },s_{n}\left( M_{\ast }\right) \right) +%
\frac{L_{\text{(\textbf{GSA}),}A_{r}}}{1-\delta -\theta _{n}}\frac{\left(
\ln n\right) ^{3}}{n}\text{ }  \notag \\
& \leq \left( \frac{1+\delta }{1-\delta }+\frac{5\theta _{n}}{\left(
1-\delta \right) ^{2}}\right) \ell \left( s_{\ast },s_{n}\left( M_{\ast
}\right) \right) +L_{\text{(\textbf{GSA}),}A_{r}}\frac{\left( \ln n\right)
^{3}}{n}\text{ .}  \label{oracle_opt_gene_proof}
\end{align}%
Inequality (\ref{oracle_opt_gene_pp}) is now proved.

\noindent It remains to prove the second part of Theorem \ref%
{theorem_opt_pen_reg_pp}. We assume that assumption (\textbf{Ap}) holds.
From Lemmas \ref{lemma_selected_model}\ and \ref{lemma_oracle_model}, we
have that for any $\frac{1}{2}>\eta >\left( 1-\beta _{+}\right) _{+}/2$ and
for all $n\geq n_{0}\left( \text{(\textbf{GSA})},C_{-},\beta _{-},\eta
,\delta \right) $, it holds on $\Omega _{n}$,%
\begin{align}
A_{\mathcal{M},+}\left( \ln n\right) ^{3}& \leq D_{%
\widehat{M}%
}\leq n^{1/2+\eta }\text{ ,}  \label{local_dim_Mhat_bis} \\
A_{\mathcal{M},+}\left( \ln n\right) ^{3}& \leq D_{M_{\ast }}\leq
n^{1/2+\eta }\text{ }.  \label{local_dim_Mstar}
\end{align}%
Now, using (\ref{lower_crit_hat_large}) and (\ref{upper_crit_star_large}),
by the same kind of computations leading to (\ref{oracle_opt_gene_proof}),
we deduce that it holds on $\Omega _{n}$, for all $n\geq n_{0}\left( \text{(%
\textbf{GSA})},C_{-},\beta _{-},\eta ,\delta \right) $,%
\begin{align*}
\ell \left( s_{\ast },s_{n}\left( 
\widehat{M}%
\right) \right) & \leq \left( \frac{1+\delta +\theta _{n}}{1-\delta -\theta
_{n}}\right) \ell \left( s_{\ast },s_{n}\left( M_{\ast }\right) \right) \\
& \leq \left( \frac{1+\delta }{1-\delta }+\frac{5\theta _{n}}{\left(
1-\delta \right) ^{2}}\right) \ell \left( s_{\ast },s_{n}\left( M_{\ast
}\right) \right) \text{ .}
\end{align*}%
Thus inequality (\ref{oracle_opt_pp}) is proved and Theorem \ref%
{theorem_opt_pen_reg_pp}\ follows. 
$\blacksquare$%

\begin{lemma}[Control on the dimension of the selected model]
\label{lemma_selected_model}Assume that (\textbf{GSA}) holds. Let $\eta \in
\left( 0,\beta _{+}/\left( 1+\beta _{+}\right) \right) $. If $n\geq
n_{0}\left( \left( \text{\textbf{GSA}}\right) ,\eta ,\delta \right) $ then,
on the event $\Omega _{n}$ defined in the proof of Theorem \ref%
{theorem_opt_pen_reg_pp}, we have%
\begin{equation}
D_{%
\widehat{M}%
}\leq n^{1/\left( 1+\beta _{+}\right) +\eta }\text{ .}
\label{control_dim_hat_gene}
\end{equation}%
If moreover (\textbf{Ap}) holds, then for all $n\geq n_{0}\left( \left( 
\text{\textbf{GSA}}\right) ,C_{-},\beta _{-},\eta ,\delta \right) $, we get
on the event $\Omega _{n}$, 
\begin{equation}
A_{\mathcal{M},+}\left( \ln n\right) ^{3}\leq D_{%
\widehat{M}%
}\leq n^{1/\left( 1+\beta _{+}\right) +\eta }\text{ .}
\label{control_dim_hat_Ap}
\end{equation}
\end{lemma}

\begin{lemma}[Control on the dimension of oracle models]
\label{lemma_oracle_model}Assume that (\textbf{GSA}) holds. Let $\eta \in
\left( 0,\beta _{+}/\left( 1+\beta _{+}\right) \right) $. If $n\geq
n_{0}\left( \left( \text{\textbf{GSA}}\right) ,\eta \right) $ then, on the
event $\Omega _{n}$ defined in the proof of Theorem \ref%
{theorem_opt_pen_reg_pp}, we have%
\begin{equation}
D_{M_{\ast }}\leq n^{1/\left( 1+\beta _{+}\right) +\eta }\text{ .}
\label{control_dim_star_gene}
\end{equation}%
If moreover (\textbf{Ap}) holds, then for all $n\geq n_{0}\left( \left( 
\text{\textbf{GSA}}\right) ,C_{-},\beta _{-},\eta \right) $, we get on the
event $\Omega _{n}$, 
\begin{equation}
A_{\mathcal{M},+}\left( \ln n\right) ^{3}\leq D_{M_{\ast }}\leq n^{1/\left(
1+\beta _{+}\right) +\eta }\text{ .}  \label{control_dim_star_Ap}
\end{equation}
\end{lemma}

\noindent \textbf{Proof of} \textbf{Lemma \ref{lemma_selected_model}}.\
Recall that $%
\widehat{M}%
$ minimizes%
\begin{equation}
\crit%
^{\prime }\left( M\right) =%
\crit%
\left( M\right) -P_{n}Ks_{\ast }=\ell \left( s_{\ast },s_{M}\right) -%
\p_{2}\left(M\right)%
+\bar{\delta}\left( M\right) +%
\pen%
\left( M\right)  \label{def_crit_lemma_selected}
\end{equation}%
over the models $M\in \mathcal{M}_{n}.$

\begin{enumerate}
\item Lower bound on $%
\crit%
^{\prime }\left( M\right) $ for small models in the case where (\textbf{Ap})
holds: let $M\in \mathcal{M}_{n}$ be such that $D_{M}<A_{\mathcal{M}%
,+}\left( \ln n\right) ^{3}.$ By (\ref{pen_id_2_pp}) and (\ref%
{def_crit_lemma_selected}), it holds%
\begin{equation*}
\crit%
^{\prime }\left( M\right) \geq \left( 1-\frac{A_{r}}{\left( \ln n\right) ^{2}%
}\right) \ell \left( s_{\ast },s_{M}\right) -%
\p_{2}\left(M\right)%
+\bar{\delta}\left( M\right) -A_{r}\frac{\left( \ln n\right) ^{3}}{n}\text{ .%
}
\end{equation*}%
We then have on $\Omega _{n}$,%
\begin{equation*}
\begin{tabular}{ll}
$\ell \left( s_{\ast },s_{M}\right) \geq C_{-}A_{\mathcal{M},+}^{-\beta
_{-}}\left( \ln n\right) ^{-3\beta _{-}}$ & $\text{by (\textbf{Ap})}$ \\ 
$%
\p_{2}\left(M\right)%
\leq L_{\text{(\textbf{GSA})}}\frac{\left( \ln n\right) ^{3}}{n}$ & $\text{%
from (\ref{p_2_small_models_1})}$ \\ 
$\bar{\delta}\left( M\right) \geq -L_{\text{(\textbf{GSA})}}\left( \sqrt{%
\frac{\ell \left( s_{\ast },s_{M}\right) \ln n}{n}}+\frac{\ln n}{n}\right) $
& $\text{from (\ref{delta_bar_small_models}).}$%
\end{tabular}%
\end{equation*}
Since by (\textbf{Ab'}), we have $0\leq \ell \left( s_{\ast },s_{M}\right)
\leq 4A^{2}$, we deduce that for all $n\geq n_{0}\left( \text{(\textbf{GSA})}%
,C_{-},\beta _{-},A_{r}\right) $,%
\begin{equation}
\crit%
^{\prime }\left( M\right) \geq \frac{C_{-}A_{\mathcal{M},+}^{-\beta _{-}}}{2}%
\left( \ln n\right) ^{-3\beta _{-}}\text{ }.  \label{lower_small_model}
\end{equation}

\item Lower bound for large models: let $M\in \mathcal{M}_{n}$ be such that $%
D_{M}\geq n^{1/\left( 1+\beta _{+}\right) +\eta }.$ From (\ref{pen_id_pp})
and (\ref{line_2}) we have on $\Omega _{n}$, for all $n\geq n_{0}\left( A_{%
\mathcal{M},+}\right) $, 
\begin{equation*}
\pen%
\left( M\right) -%
\p_{2}\left(M\right)%
\geq \mathbb{E}\left[ 
\p_{2}\left(M\right)%
\right] -\left( \delta +L_{\text{(\textbf{GSA})}}\varepsilon _{n}^{2}\left(
M\right) \right) \left( \ell \left( s_{\ast },s_{M}\right) +\mathbb{E}\left[ 
\p_{2}\left(M\right)%
\right] \right) \text{ .}
\end{equation*}%
Using (\textbf{P2}) and the fact that $D_{M}\geq n^{1/\left( 1+\beta
_{+}\right) +\eta }$ in (\ref{def_epsilon_selection}), we deduce that for
all $n\geq n_{0}\left( \text{(\textbf{GSA})},\eta ,\delta ,\beta _{+}\right) 
$, $L_{\text{(\textbf{GSA})}}\varepsilon _{n}^{2}\left( M\right) \leq \frac{1%
}{2}\left( 1-\delta \right) $ and as by (\textbf{An}), $\mathcal{K}%
_{1,M}\geq 2\sigma _{\min }$, we also deduce from Lemma \ref{mean_emp_risk}
that for all $n\geq n_{0}\left( \text{(\textbf{GSA})},\eta \right) $, $%
\mathbb{E}\left[ 
\p_{2}\left(M\right)%
\right] \geq \frac{\sigma _{\min }^{2}}{2}\frac{D_{M}}{n}$. Consequently, it
holds for all $n\geq n_{0}\left( \text{(\textbf{GSA})},\eta ,\delta ,\beta
_{+}\right) $, 
\begin{equation}
\pen%
\left( M\right) -%
\p_{2}\left(M\right)%
\geq \frac{\sigma _{\min }^{2}}{4}\left( 1-\delta \right) \frac{D_{M}}{n}%
-C_{+}D_{M}^{-\beta _{+}}\geq \left( 1-\delta \right) L_{\text{(\textbf{GSA})%
}}n^{-\frac{\beta _{+}}{1+\beta _{+}}+\eta }
\label{lower_pen_large_selected}
\end{equation}%
From (\ref{line_4}) it holds on $\Omega _{n}$,%
\begin{equation}
\bar{\delta}\left( M\right) \geq -L_{\text{(\textbf{GSA})}}\left( \sqrt{%
\frac{\ell \left( s_{\ast },s_{M}\right) \ln n}{n}}+\frac{\ln n}{n}\right)
\geq -L_{\text{(\textbf{GSA})}}\left( n^{-\frac{1+2\beta _{+}}{2\left(
1+\beta _{+}\right) }}\sqrt{\ln n}+\frac{\ln n}{n}\right) \text{ .}
\label{lower_delta_large_selected}
\end{equation}%
Hence, we deduce from (\ref{def_crit_lemma_selected}), (\ref%
{lower_pen_large_selected}) and (\ref{lower_delta_large_selected}) that we
have on $\Omega _{n}$, for all $n\geq n_{0}\left( \text{(\textbf{GSA})},\eta
,\delta ,\beta _{+}\right) $, 
\begin{equation}
\crit%
^{\prime }\left( M\right) \geq \left( 1-\delta \right) L_{\text{(\textbf{GSA}%
)}}n^{-\frac{\beta _{+}}{1+\beta _{+}}+\eta }\text{ }.
\label{lower_large_model}
\end{equation}

\item A better model exists for $%
\crit%
^{\prime }\left( M\right) $: from (\textbf{P3}), there exists $M_{0}\in 
\mathcal{M}_{n}$ such that $n^{1/\left( 1+\beta _{+}\right) }\leq
D_{M_{0}}\leq c_{rich}n^{1/\left( 1+\beta _{+}\right) }.$ Then, for all $%
n\geq n_{0}\left( \text{(\textbf{GSA})},\eta \right) $,%
\begin{equation*}
A_{\mathcal{M},+}\left( \ln n\right) ^{3}\leq n^{1/\left( 1+\beta
_{+}\right) }\leq D_{M_{0}}\leq c_{rich}n^{1/\left( 1+\beta _{+}\right)
}\leq n^{1/\left( 1+\beta _{+}\right) +\eta }\text{ }.
\end{equation*}%
Using (\textbf{Ap}$_{u}$),%
\begin{equation}
\ell \left( s_{\ast },s_{M_{0}}\right) \leq C_{+}n^{-\beta _{+}/\left(
1+\beta _{+}\right) }\text{ }.  \label{upper_bias_M0_selected}
\end{equation}%
By (\ref{line_3_ter}), we have on $\Omega _{n}$, for all $n\geq n_{0}\left( 
\text{(\textbf{GSA})},\eta \right) $,%
\begin{equation}
\left\vert \bar{\delta}\left( M_{0}\right) \right\vert \leq \frac{\ell
\left( s_{\ast },s_{M_{0}}\right) }{\sqrt{D_{M_{0}}}}+L_{\text{(\textbf{GSA})%
}}\frac{\ln n}{\sqrt{D_{M_{0}}}}\mathbb{E}\left[ \text{p}_{2}\left(
M_{0}\right) \right] \leq L_{\text{(\textbf{GSA})}}n^{-\frac{1+2\beta _{+}}{%
2\left( 1+\beta _{+}\right) }}\ln \left( n\right)
\label{upper_delta_M0_selected}
\end{equation}%
and by (\ref{pen_id_pp}),%
\begin{equation*}
\pen%
\left( M_{0}\right) \leq 3\left( \ell \left( s_{\ast },s_{M_{0}}\right) +%
\mathbb{E}\left[ \text{p}_{2}\left( M_{0}\right) \right] \right) \leq L_{%
\text{(\textbf{GSA})}}n^{-\beta _{+}/\left( 1+\beta _{+}\right) }\text{ .}
\end{equation*}%
Consequently, we have on $\Omega _{n}$, for all $n\geq n_{0}\left( \text{(%
\textbf{GSA})},\eta \right) $,%
\begin{align}
\crit%
^{\prime }\left( M_{0}\right) & \leq \ell \left( s_{\ast },s_{M_{0}}\right)
+\left\vert \bar{\delta}\left( M_{0}\right) \right\vert +%
\pen%
\left( M_{0}\right)  \notag \\
& \leq L_{\text{(\textbf{GSA})}}n^{-\beta _{+}/\left( 1+\beta _{+}\right) }%
\text{ .}  \label{upper_bon_modele}
\end{align}
\end{enumerate}

\noindent To conclude, notice that the upper bound (\ref{upper_bon_modele})
is smaller than the lower bound given in (\ref{lower_large_model}) for all $%
n\geq n_{0}\left( \text{(\textbf{GSA})},\eta ,\delta \right) $. Hence,
points 2 and 3\textbf{\ }above\textbf{\ }yield inequality (\ref%
{control_dim_hat_gene}). Moreover, the upper bound (\ref{upper_bon_modele})
is smaller than lower bounds given in (\ref{lower_small_model}), derived by
using (\textbf{Ap}), and (\ref{lower_large_model}), for all $n\geq
n_{0}\left( \text{(\textbf{GSA})},C_{-},\beta _{-},\eta ,\delta \right) $.
This thus gives (\ref{control_dim_hat_Ap}) and Lemma \ref%
{lemma_selected_model}\ is proved. 
$\blacksquare$%

\noindent \textbf{Proof of} \textbf{Lemma \ref{lemma_oracle_model}}. By
definition, $M_{\ast }$ minimizes 
\begin{equation*}
\ell \left( s_{\ast },s_{n}\left( M\right) \right) =\ell \left( s_{\ast
},s_{M}\right) +%
\p_{1}\left(M\right)%
\end{equation*}%
over the models $M\in \mathcal{M}_{n}.$

\begin{enumerate}
\item Lower bound on $\ell \left( s_{\ast },s_{n}\left( M\right) \right) $
for small models: let $M\in \mathcal{M}_{n}$ be such that $D_{M}<A_{\mathcal{%
M},+}\left( \ln n\right) ^{3}.$ In this case we have%
\begin{equation}
\ell \left( s_{\ast },s_{n}\left( M\right) \right) \geq \ell \left( s_{\ast
},s_{M}\right) \geq C_{-}A_{\mathcal{M},+}^{-\beta _{-}}\left( \ln n\right)
^{-3\beta _{-}}\text{ by (\textbf{Ap}).}  \label{lower_small_oracle}
\end{equation}

\item Lower bound of $\ell \left( s_{\ast },s_{n}\left( M\right) \right) $
for large models: let $M\in \mathcal{M}_{n}$ be such that $D_{M}\geq
n^{1/\left( 1+\beta _{+}\right) +\eta }.$ From (\ref{line_1}) we get on $%
\Omega _{n}$,%
\begin{equation*}
\p_{1}\left(M\right)%
\geq \left( 1-L_{\text{(\textbf{GSA})}}\varepsilon _{n}\left( M\right)
\right) \mathbb{E}\left[ \text{p}_{2}\left( M\right) \right] \text{ .}
\end{equation*}%
Using (\textbf{P2}) and the fact that $D_{M}\geq n^{1/\left( 1+\beta
_{+}\right) +\eta }$ in (\ref{def_epsilon_selection}), we deduce that for
all $n\geq n_{0}\left( \text{(\textbf{GSA})},\eta \right) $, $L_{\text{(%
\textbf{GSA})}}\varepsilon _{n}\left( M\right) \leq \frac{1}{2}$ and as by (%
\textbf{An}), $\mathcal{K}_{1,M}\geq 2\sigma _{\min }$ we also deduce from
Lemma \ref{mean_emp_risk} that for all $n\geq n_{0}\left( \text{(\textbf{GSA}%
)},\eta \right) $, $\mathbb{E}\left[ \text{p}_{2}\left( M\right) \right]
\geq \frac{\sigma _{\min }^{2}}{2}\frac{D_{M}}{n}$. Consequently, it holds
for all $n\geq n_{0}\left( \text{(\textbf{GSA})},\eta \right) $, on the
event $\Omega _{n}$, 
\begin{equation}
\ell \left( s_{\ast },s_{n}\left( M\right) \right) \geq 
\p_{1}\left(M\right)%
\geq \frac{\sigma _{\min }^{2}}{4}\frac{D_{M}}{n}\geq \frac{\sigma _{\min
}^{2}}{4}n^{-\beta _{+}/\left( 1+\beta _{+}\right) +\eta }\text{ }.
\label{lower_large_oracle}
\end{equation}

\item A better model exists for $\ell \left( s_{\ast },s_{n}\left( M\right)
\right) $: from (\textbf{P3}), there exists $M_{0}\in \mathcal{M}_{n}$ such
that $n^{1/\left( 1+\beta _{+}\right) }\leq D_{M_{0}}\leq
c_{rich}n^{1/\left( 1+\beta _{+}\right) }.$ Moreover, for all $n\geq
n_{0}\left( \text{(\textbf{GSA})},\eta \right) $,%
\begin{equation*}
A_{\mathcal{M},+}\left( \ln n\right) ^{3}\leq n^{1/\left( 1+\beta
_{+}\right) }\leq D_{M_{0}}\leq c_{rich}n^{1/\left( 1+\beta _{+}\right)
}\leq n^{1/\left( 1+\beta _{+}\right) +\eta }\text{ }.
\end{equation*}%
Using (\textbf{Ap}$_{u}$),%
\begin{equation*}
\ell \left( s_{\ast },s_{M_{0}}\right) \leq C_{+}n^{-\beta _{+}/\left(
1+\beta _{+}\right) }
\end{equation*}%
and by (\ref{line_1})%
\begin{equation*}
\text{p}_{1}\left( M_{0}\right) \leq \left( 1+L_{\text{(\textbf{GSA})}%
}\varepsilon _{n}\left( M\right) \right) \mathbb{E}\left[ \text{p}_{2}\left(
M_{0}\right) \right] \text{ .}
\end{equation*}%
Hence, as $\mathcal{K}_{1,M}\leq 6A$ by (\textbf{Ab'}) and as, by (\ref%
{def_epsilon_selection}), for all $n\geq n_{0}\left( \text{\textbf{GSA}}%
\right) $ it holds $\varepsilon _{n}\left( M\right) \leq 1$, we deduce from
Lemma \ref{mean_emp_risk} that for all $n\geq n_{0}\left( \text{\textbf{GSA}}%
\right) $, on the event $\Omega _{n}$,%
\begin{equation*}
\text{p}_{1}\left( M_{0}\right) \leq L_{\text{(\textbf{GSA})}}\frac{D_{M}}{n}%
\leq L_{\text{(\textbf{GSA})}}n^{-\beta _{+}/\left( 1+\beta _{+}\right) }%
\text{ }.
\end{equation*}%
Consequently, on $\Omega _{n}$, for all $n\geq n_{0}\left( \text{(\textbf{GSA%
})},\eta \right) $,%
\begin{align}
\ell \left( s_{\ast },s_{n}\left( M_{0}\right) \right) & =\ell \left(
s_{\ast },s_{M_{0}}\right) +\text{p}_{1}\left( M_{0}\right)  \notag \\
& \leq L_{\text{(\textbf{GSA})}}n^{-\beta _{+}/\left( 1+\beta _{+}\right) }%
\text{ }.  \label{upper_reasonable_model}
\end{align}
\end{enumerate}

\noindent The upper bound (\ref{upper_reasonable_model}) is smaller than the
lower bound (\ref{lower_large_oracle}) for all $n\geq n_{0}\left( \text{(%
\textbf{GSA})},\eta \right) $, and this gives (\ref{control_dim_star_gene}).
If (\textbf{Ap}) holds, then the upper bound (\ref{upper_reasonable_model})
is smaller than the lower bounds (\ref{lower_small_oracle}) and (\ref%
{lower_large_oracle}) for all $n\geq n_{0}\left( \text{(\textbf{GSA})}%
,C_{-},\beta _{-},\eta \right) $, which proves (\ref{control_dim_star_Ap})
and allows to conclude the proof of Lemma \ref{lemma_oracle_model}.\ 
$\blacksquare$%

\subsection*{}%
\textbf{Proof of Theorem \ref{theorem_min_pen_reg_pp}. }As in the proof of
Theorem \ref{theorem_opt_pen_reg_pp}, we consider the event $\Omega
_{n}^{\prime }$ of probability at least $1-L_{c_{\mathcal{M}},A_{p}}n^{-2}$
for all $n\geq n_{0}\left( \text{\textbf{GSA}}\right) $, on which: (\ref%
{majo_pen_pp}) holds and

\begin{itemize}
\item For all models $M\in \mathcal{M}_{n}$ of dimension $D_{M}$ such that $%
A_{\mathcal{M},+}\left( \ln n\right) ^{2}\leq D_{M}$, 
\begin{align}
\left\vert 
\p_{1}\left(M\right)%
-\mathbb{E}\left[ \text{p}_{2}\left( M\right) \right] \right\vert & \leq L_{%
\text{(\textbf{GSA})}}\varepsilon _{n}\left( M\right) \mathbb{E}\left[ \text{%
p}_{2}\left( M\right) \right] \text{ ,}  \label{line_1_min} \\
\left\vert \text{p}_{2}\left( M\right) -\mathbb{E}\left[ \text{p}_{2}\left(
M\right) \right] \right\vert & \leq L_{\text{(\textbf{GSA})}}\varepsilon
_{n}^{2}\left( M\right) \mathbb{E}\left[ \text{p}_{2}\left( M\right) \right] 
\text{ .}  \label{line_2_min_bis}
\end{align}

\item For all models $M\in \mathcal{M}_{n}$ with $D_{M}\leq A_{\mathcal{M}%
,+}\left( \ln n\right) ^{2}$,%
\begin{equation}
\text{p}_{2}\left( M\right) \leq L_{\text{(\textbf{GSA})}}\frac{\left( \ln
n\right) ^{2}}{n}\text{ .}  \label{p_2_small_models_min}
\end{equation}

\item For every $M\in \mathcal{M}_{n}$,%
\begin{equation}
\left\vert \bar{\delta}\left( M\right) \right\vert \leq L_{\text{(\textbf{GSA%
})}}\left( \sqrt{\frac{\ell \left( s_{\ast },s_{M}\right) \ln n}{n}}+\frac{%
\ln n}{n}\right) \text{ .}  \label{line_4_min}
\end{equation}
\end{itemize}

\noindent Let $d\in \left( 0,1\right) $ to be chosen later.

\noindent \textbf{Lower bound on} $D_{%
\widehat{M}%
}$.\ Let us recall that $%
\widehat{M}%
$ minimizes%
\begin{equation}
\crit%
^{\prime }\left( M\right) =%
\crit%
\left( M\right) -P_{n}Ks_{\ast }=\ell \left( s_{\ast },s_{M}\right) -%
\p_{2}\left(M\right)%
+\bar{\delta}\left( M\right) +%
\pen%
\left( M\right) \text{ }.  \label{def_crit_prime}
\end{equation}

\begin{enumerate}
\item Lower bound on $%
\crit%
^{\prime }\left( M\right) $ for 
``%
small%
''
models: assume that $M\in \mathcal{M}_{n}$ and 
\begin{equation*}
D_{M}\leq dA_{rich}n\left( \ln n\right) ^{-2}\text{ }.
\end{equation*}%
We have 
\begin{equation}
\ell \left( s_{\ast },s_{M}\right) +%
\pen%
\left( M\right) \geq 0  \label{crit_part1}
\end{equation}%
and from (\ref{line_4_min}), as $\ell \left( s_{\ast },s_{M}\right) \leq
4A^{2}$ by (\textbf{Ab'}), we get on $\Omega _{n}^{\prime }$, for all $n\geq
n_{0}\left( \text{(\textbf{GSA}),}d\right) $,%
\begin{align}
\bar{\delta}\left( M\right) & \geq -L_{\text{(\textbf{GSA})}}\left( \sqrt{%
\frac{\ell \left( s_{\ast },s_{M}\right) \ln n}{n}}+\frac{\ln n}{n}\right) 
\text{ }  \notag \\
& \geq -L_{\text{(\textbf{GSA})}}\sqrt{\frac{\ln n}{n}}  \notag \\
& \geq -d\times A^{2}A_{rich}\left( \ln n\right) ^{-2}\text{ .}
\label{crit_part2}
\end{align}%
Then, if $D_{M}\geq A_{\mathcal{M},+}\left( \ln n\right) ^{2}$, as $\mathcal{%
K}_{1,M}\leq 6A$ by (\textbf{Ab'}) and as, by (\ref{def_epsilon_selection}),
for all $n\geq n_{0}\left( \text{\textbf{GSA}}\right) $ it holds $L_{\text{(%
\textbf{GSA})}}\varepsilon _{n}\left( M\right) \leq 1$, we deduce from (\ref%
{line_2_min_bis}) and Lemma \ref{mean_emp_risk} that for all $n\geq
n_{0}\left( \text{(\textbf{GSA}),}d\right) $, 
\begin{equation*}
\p_{2}\left(M\right)%
\leq 2\mathbb{E}\left[ 
\p_{2}\left(M\right)%
\right] \leq 36A^{2}\frac{D_{M}}{n}\leq d\times 36A^{2}A_{rich}\left( \ln
n\right) ^{-2}\text{ }.
\end{equation*}%
Whenever $D_{M}\leq A_{\mathcal{M},+}\left( \ln n\right) ^{2}$, (\ref%
{p_2_small_models_min}) gives that, for all $n\geq n_{0}\left( \text{(%
\textbf{GSA}),}d\right) $, on the event $\Omega _{n}^{\prime }$, 
\begin{equation*}
\p_{2}\left(M\right)%
\leq L_{\text{(\textbf{GSA})}}\frac{\left( \ln n\right) ^{2}}{n}\leq d\times
36A^{2}A_{rich}\left( \ln n\right) ^{-2}\text{ }.
\end{equation*}%
Hence, we have checked that for all $n\geq n_{0}\left( \text{(\textbf{GSA})}%
,d\right) $, on the event $\Omega _{n}^{\prime }$,%
\begin{equation}
-%
\p_{2}\left(M\right)%
\geq -d\times 36A^{2}A_{rich}\left( \ln n\right) ^{-2}\text{ },
\label{crit_part3}
\end{equation}%
and finally, by using (\ref{crit_part1}), (\ref{crit_part2}) and (\ref%
{crit_part3}) in (\ref{def_crit_prime}), we deduce that on $\Omega
_{n}^{\prime }$, for all $n\geq n_{0}\left( \text{(\textbf{GSA})},d\right) $,%
\begin{equation}
\crit%
^{\prime }\left( M\right) \geq -d\times 37A^{2}A_{rich}\left( \ln n\right)
^{-2}\text{ .}  \label{lower_crit_small}
\end{equation}

\item There exists a better model for $%
\crit%
^{\prime }\left( M\right) $. By (\textbf{P3}), for all $n\geq n_{0}\left( A_{%
\mathcal{M},+},A_{rich}\right) $ a model $M_{1}\in \mathcal{M}_{n}$ exists
such that%
\begin{equation*}
A_{\mathcal{M},+}\left( \ln n\right) ^{2}\leq \frac{A_{rich}n}{\left( \ln
n\right) ^{2}}\leq D_{M_{1}}\text{ }.
\end{equation*}%
We then have on $\Omega _{n}^{\prime }$,%
\begin{equation*}
\begin{tabular}{ll}
$\ell \left( s_{\ast },s_{M_{1}}\right) \leq A_{rich}^{-\beta _{+}}\left(
\ln n\right) ^{2\beta _{+}}n^{-\beta _{+}}$ & $\text{by (\textbf{Ap}}_{u}%
\text{)}$ \\ 
$\text{p}_{2}\left( M_{1}\right) \geq \left( 1-L_{\text{(\textbf{GSA})}%
}\varepsilon _{n}^{2}\left( M_{1}\right) \right) \mathbb{E}\left[ \text{p}%
_{2}\left( M_{1}\right) \right] $ & $\text{by (\ref{line_2_min_bis})}$ \\ 
$%
\pen%
\left( M_{1}\right) \leq A_{%
\pen%
}\mathbb{E}\left[ \text{p}_{2}\left( M_{1}\right) \right] $ & $\text{by (\ref%
{majo_pen_pp})}$ \\ 
$\left\vert \bar{\delta}\left( M_{1}\right) \right\vert \leq L_{\text{(%
\textbf{GSA})}}\sqrt{\ln \left( n\right) /n}$ & $\text{by (\ref{line_4_min})
and (\textbf{Ab'})}$%
\end{tabular}%
\end{equation*}%
and therefore, 
\begin{equation}
\crit%
^{\prime }\left( M_{1}\right) \leq \left( -1+A_{%
\pen%
}+L_{\text{(\textbf{GSA})}}\varepsilon _{n}^{2}\left( M_{1}\right) \right) 
\mathbb{E}\left[ \text{p}_{2}\left( M_{1}\right) \right] +L_{\text{(\textbf{%
GSA})}}\sqrt{\frac{\ln n}{n}}+A_{rich}^{-\beta _{+}}\frac{\left( \ln
n\right) ^{2\beta _{+}}}{n^{\beta _{+}}}\text{ }.  \label{upper_crit_1}
\end{equation}%
Hence, as $-1+A_{%
\pen%
}<0$, and as by (\ref{def_epsilon_selection}), (\textbf{An}) and Lemma \ref%
{mean_emp_risk} it holds for all $n\geq n_{0}\left( \text{(\textbf{GSA})},A_{%
\pen%
}\right) $%
\begin{equation*}
L_{\text{(\textbf{GSA})}}\varepsilon _{n}^{2}\left( M_{1}\right) \leq \frac{%
1-A_{%
\pen%
}}{2}\text{ \ \ \ and \ \ \ }\mathbb{E}\left[ \text{p}_{2}\left(
M_{1}\right) \right] \geq \frac{\sigma _{\min }^{2}}{2}\frac{D_{M}}{n}\geq 
\frac{\sigma _{\min }^{2}A_{rich}}{2}\left( \ln n\right) ^{-2}\text{ ,}
\end{equation*}%
we deduce from (\ref{upper_crit_1}) that on $\Omega _{n}^{\prime }$, for all 
$n\geq n_{0}\left( \text{(\textbf{GSA})},A_{%
\pen%
}\right) $,%
\begin{equation}
\crit%
^{\prime }\left( M_{1}\right) \leq -\frac{1}{4}\left( 1-A_{%
\pen%
}\right) \sigma _{\min }^{2}A_{rich}\left( \ln n\right) ^{-2}\text{ }.
\label{upper_crit_2}
\end{equation}
\end{enumerate}

\noindent Now, by taking 
\begin{equation}
0<d=\left( \frac{1}{149}\left( 1-A_{%
\pen%
}\right) \left( \frac{\sigma _{\min }}{A}\right) ^{2}\right) \wedge \frac{1}{%
2}<1  \label{value_d_reg}
\end{equation}%
and by comparing (\ref{lower_crit_small}) and (\ref{upper_crit_2}), we
deduce that on $\Omega _{n}^{\prime }$, for all $n\geq n_{0}\left( \text{(%
\textbf{GSA})},A_{%
\pen%
}\right) $, for all $M\in \mathcal{M}_{n}$ such that $D_{M}\leq
dA_{rich}n\left( \ln n\right) ^{-2}$, 
\begin{equation*}
\crit%
^{\prime }\left( M_{1}\right) <%
\crit%
^{\prime }\left( M\right)
\end{equation*}%
and so 
\begin{equation}
D_{%
\widehat{M}%
}>dA_{rich}n\left( \ln n\right) ^{-2}\text{ }.  \label{lower_dim_hat}
\end{equation}

\noindent \textbf{Excess} \textbf{Loss of} $s_{n}\left( 
\widehat{M}%
\right) $.\ We take $d$ with the value given in (\ref{value_d_reg}). First
notice that for all $n\geq n_{0}\left( A_{\mathcal{M},+},A_{rich},d\right) ,$
we have $dA_{rich}n\left( \ln n\right) ^{-2}\geq A_{\mathcal{M},+}\left( \ln
n\right) ^{2}$. Hence, for all $M\in \mathcal{M}_{n}$ such that $D_{M}\geq
dA_{rich}n\left( \ln n\right) ^{-2}$, by (\ref{def_epsilon_selection}), (%
\textbf{P2}), (\textbf{An}) and Lemma \ref{mean_emp_risk}, it holds on $%
\Omega _{n}^{\prime }$ for all $n\geq n_{0}\left( \text{(\textbf{GSA})},A_{%
\pen%
}\right) $, using (\ref{line_1_min}),%
\begin{equation*}
\ell \left( s_{\ast },s_{n}\left( M\right) \right) \geq 
\p_{1}\left(M\right)%
\geq \frac{\sigma _{\min }^{2}}{2}\frac{D_{M}}{n}\geq \frac{d\sigma _{\min
}^{2}A_{rich}}{2}\left( \ln n\right) ^{-2}\text{ }.
\end{equation*}%
By (\ref{lower_dim_hat}), we thus get that on $\Omega _{n}^{\prime }$, for
all $n\geq n_{0}\left( \text{(\textbf{GSA})},A_{%
\pen%
}\right) $, 
\begin{equation}
\ell \left( s_{\ast },s_{n}\left( 
\widehat{M}%
\right) \right) \geq \frac{d\sigma _{\min }^{2}A_{rich}}{2}\left( \ln
n\right) ^{-2}\text{ }.  \label{excess_risk_large}
\end{equation}%
Moreover, the model $M_{0}$ defined in (\textbf{P3}) satisfies, for all $%
n\geq n_{0}\left( \text{\textbf{GSA}}\right) $,%
\begin{equation*}
A_{\mathcal{M},+}\left( \ln n\right) ^{3}\leq n^{1/\left( 1+\beta
_{+}\right) }\leq D_{M_{0}}\leq c_{rich}n^{1/\left( 1+\beta _{+}\right) }
\end{equation*}%
and so using (\textbf{Ap}$_{u}$),%
\begin{equation*}
\ell \left( s_{\ast },s_{M_{0}}\right) \leq C_{+}n^{-\beta _{+}/\left(
1+\beta _{+}\right) }\text{ .}
\end{equation*}%
In addition, by (\ref{line_1}),%
\begin{equation*}
\p_{1}\left(M\right)%
\leq \left( 1+L_{\text{(\textbf{GSA})}}\varepsilon _{n}\left( M\right)
\right) \mathbb{E}\left[ \text{p}_{2}\left( M\right) \right] \text{ .}
\end{equation*}%
Hence, as $\mathcal{K}_{1,M}\leq 6A$ by (\textbf{Ab'}) and as, by (\ref%
{def_epsilon_selection}), for all $n\geq n_{0}\left( \text{\textbf{GSA}}%
\right) $ it holds $\varepsilon _{n}\left( M\right) \leq 1$, we deduce from
Lemma \ref{mean_emp_risk} that for all $n\geq n_{0}\left( \text{\textbf{GSA}}%
\right) $,%
\begin{equation*}
\p_{1}\left(M\right)%
\leq L_{\text{(\textbf{GSA})}}\frac{D_{M}}{n}\leq L_{\text{(\textbf{GSA})}%
}n^{-\beta _{+}/\left( 1+\beta _{+}\right) }\text{ }.
\end{equation*}%
Consequently, for all $n\geq n_{0}\left( \text{\textbf{GSA}}\right) $,%
\begin{equation}
\ell \left( s_{\ast },s_{n}\left( M_{0}\right) \right) \leq L_{\text{(%
\textbf{GSA})}}n^{-\beta _{+}/\left( 1+\beta _{+}\right) }
\label{excess_risk_M0}
\end{equation}%
and the ratio between the two bounds (\ref{excess_risk_large}) and (\ref%
{excess_risk_M0}) is larger than $n^{\beta _{+}/\left( 1+\beta _{+}\right)
}\left( \ln n\right) ^{-3}$ for all $n\geq n_{0}\left( L_{\text{(\textbf{GSA}%
)}},A_{%
\pen%
}\right) $, which yields (\ref{bad_oracle_min_pen_pp}).\ 
$\blacksquare$%

\subsection{Proofs related to Section \protect\ref{section_hold_out}\label%
{section_proof_hold_out}}

Theorem \ref{theorem_pen_n1n2_pp} is a straightforward consequence of the
following result, that will be proved below.

\begin{theorem}
\label{theorem_pen_n1n2}Assume that (\textbf{GSA}) holds. With the notations
of Section \ref{section_hold_out}, assume moreover that there exist $c\in
\left( 0,1\right) $ such that $nc\leq n_{1}<n$ and $\tau \in \left(
1,3\right) $ satisfying $n\left( \ln n\right) ^{\tau }/D_{M}\leq n_{2}\leq
n\left( 1-c\right) $ for all $M\in \mathcal{M}_{n}$ such that $A_{\mathcal{M}%
,+}\left( \ln n\right) ^{3}\leq D_{M}\leq A_{\mathcal{M},+}n/\left( \ln
n\right) ^{2}$. Take $n_{2}=n\left( 1-c\right) $\ if $D_{M}\leq A_{\mathcal{M%
},+}\left( \ln n\right) ^{3}$. Define for all $M\in \mathcal{M}_{n}$, 
\begin{equation*}
\pen%
_{ho}\left( M\right) =\frac{n_{1}}{n}\left( P_{n_{2}}\left( Ks_{n_{1}}\left(
M\right) \right) -P_{n_{1}}\left( Ks_{n_{1}}\left( M\right) \right) \right) 
\text{ .}
\end{equation*}%
Then, for any $\eta \in \left( 0,\beta _{+}/\left( 1+\beta _{+}\right)
\right) $, there exist an integer $n_{0}$ depending on $c,\eta $ and on
constants in (\textbf{GSA)}, a positive constant $A_{6}$ only depending on $%
c_{\mathcal{M}}$ given in (\textbf{GSA}), two positive constants $A_{7}$ and 
$A_{8}$ only depending on constants in (\textbf{GSA}) and a sequence 
\begin{equation*}
\theta _{n}\leq \frac{A_{7}}{\left( \ln n\right) ^{1/4}\wedge \left( \ln
n\right) ^{\left( \tau -1\right) /2}}
\end{equation*}%
such that it holds for all $n\geq n_{0}\left( \left( \text{\textbf{GSA}}%
\right) ,c,\eta \right) $, with probability at least $1-A_{6}n^{-2}$,%
\begin{equation*}
D_{%
\widehat{M}%
_{n_{1}}}\leq n^{\eta +1/\left( 1+\beta _{+}\right) }
\end{equation*}%
and%
\begin{equation}
\ell \left( s_{\ast },s_{n}\left( 
\widehat{M}%
_{n_{1}}\right) \right) \leq \left( 1+\theta _{n}\right) \ell \left( s_{\ast
},s_{n}\left( M_{\ast }\right) \right) +A_{8}\frac{\left( \ln n\right) ^{3}}{%
n}\text{ .}  \label{oracle_opt_gene_n1}
\end{equation}%
Assume that in addition (\textbf{Ap}) holds (see Theorem \ref%
{theorem_opt_pen_reg_pp}). Then it holds for all $n\geq n_{0}\left( \left( 
\text{\textbf{GSA}}\right) ,C_{-},\beta _{-},\eta ,c\right) $, with
probability at least $1-A_{6}n^{-2}$,%
\begin{equation*}
A_{\mathcal{M},+}\left( \ln n\right) ^{3}\leq D_{%
\widehat{M}%
_{n_{1}}}\leq n^{\eta +1/\left( 1+\beta _{+}\right) }
\end{equation*}%
and%
\begin{equation}
\ell \left( s_{\ast },s_{n}\left( 
\widehat{M}%
_{n_{1}}\right) \right) \leq \left( 1+\theta _{n}\right) \inf_{M\in \mathcal{%
M}_{n}}\left\{ \ell \left( s_{\ast },s_{n}\left( M\right) \right) \right\} 
\text{ .}  \label{oracle_opt_n1}
\end{equation}
\end{theorem}

\begin{lemma}
\label{Lemma_p1n2}Assume that (\textbf{GSA}) holds. Let $c\in \left(
0,1\right) $, $\tau \in \left( 1,3\right) $ and $\left( n_{1},n_{2}\right)
\in \mathbb{N}_{\ast }^{2}$. We assume that $nc\leq n_{1}<n$ and set $%
n_{2}=n-n_{1}$. Then there exists $L=L_{\left( \text{\textbf{GSA}}\right) 
\text{,}c}>0$ such that for all $M\in \mathcal{M}_{n}$ satisfying $D_{M}\geq
A_{\mathcal{M},+}\left( \ln n\right) ^{2}$, for all $n\geq n_{0}\left(
\left( \text{\textbf{GSA}}\right) \text{,}c\right) $, it holds%
\begin{gather}
\mathbb{P}\left( \left\vert P_{n_{2}}\left( Ks_{n_{1}}\left( M\right)
-Ks_{M}\right) -P\left( Ks_{n_{1}}\left( M\right) -Ks_{M}\right) \right\vert
\geq L\frac{\sqrt{\left( D_{M}\vee \ln n\right) \left( \ln n\right) \left(
\left( \ln n\right) \left( \ln n_{1}\right) +n_{2}\right) }}{n_{2}\sqrt{n_{1}%
}}\right)  \notag \\
\leq 12n^{-2-\alpha _{\mathcal{M}}}\text{ .}  \label{p1n2_1}
\end{gather}%
Now, let us assume that $n\left( \ln n\right) ^{\tau }/D_{M}\leq n_{2}\leq
n\left( 1-c\right) $ if $A_{\mathcal{M},+}\left( \ln n\right) ^{3}\leq
D_{M}\leq A_{\mathcal{M},+}n/\left( \ln n\right) ^{2}$ and $n_{2}=n\left(
1-c\right) $\ if $D_{M}\leq A_{\mathcal{M},+}\left( \ln n\right) ^{3}$. If $%
A_{\mathcal{M},+}\left( \ln n\right) ^{3}\leq D_{M}$, then by setting%
\begin{equation}
\varepsilon _{n}^{1,2}\left( M\right) =L\frac{n\sqrt{\ln n\left( \left( \ln
n\right) \left( \ln n_{1}\right) +n_{2}\right) }}{n_{2}\sqrt{n_{1}D_{M}}}%
\leq \frac{L}{\left( \ln n\right) ^{\left( \tau -1\right) /2}}\text{ ,}
\label{def_epsilon_1_2}
\end{equation}%
we have for all $n\geq n_{0}\left( \left( \text{\textbf{GSA}}\right)
,c\right) $,%
\begin{equation}
\mathbb{P}\left( \left\vert P_{n_{2}}\left( Ks_{n_{1}}\left( M\right)
-Ks_{M}\right) -P\left( Ks_{n_{1}}\left( M\right) -Ks_{M}\right) \right\vert
\geq \varepsilon _{n}^{1,2}\left( M\right) \mathbb{E}\left[ 
\p_{2}\left(M\right)%
\right] \right) \leq 12n^{-2-\alpha _{\mathcal{M}}}\text{ .}  \label{p1n2_3}
\end{equation}%
If $D_{M}\leq A_{\mathcal{M},+}\left( \ln n\right) ^{3}$, we obtain%
\begin{equation}
\mathbb{P}\left( \left\vert P_{n_{2}}\left( Ks_{n_{1}}\left( M\right)
-Ks_{M}\right) -P\left( Ks_{n_{1}}\left( M\right) -Ks_{M}\right) \right\vert
\geq L\frac{\left( \ln n\right) ^{2}}{n}\right) \leq 12n^{-2-\alpha _{%
\mathcal{M}}}\text{ .}  \label{p1n2_4}
\end{equation}
\end{lemma}

\begin{proof}
By Bernstein's inequality (see Corollary 2.10 in \cite{Massart:07}) applied
to the sum of $\left( s_{n_{1}}\left( M\right) \right) \left( \xi
_{i}\right) $ conditionally to $\left( \xi _{j}\right) _{j\in I_{1}}$, we
get that for all $x>0$, it holds%
\begin{equation}
\mathbb{P}\left( \left\vert P_{n_{2}}\left( Ks_{n_{1}}\left( M\right)
-Ks_{M}\right) -P\left( Ks_{n_{1}}\left( M\right) -Ks_{M}\right) \right\vert
\geq x\left\vert \left( \xi _{j}\right) ,\text{ }j\in I_{1}\right. \right)
\leq 2\exp \left( -\frac{nx^{2}}{2\left( v_{1}+b_{1}x/3\right) }\right) 
\text{ ,}  \label{condi_ber}
\end{equation}%
where 
\begin{equation*}
v_{1}=\mathbb{E}_{\xi }\left[ \left( Ks_{n_{1}}\left( M\right) \left( \xi
\right) -Ks_{M}\left( \xi \right) \right) ^{2}\right]
\end{equation*}%
and $b_{1}=\left\Vert Ks_{n_{1}}\left( M\right) -Ks_{M}\right\Vert _{\infty
} $. We have%
\begin{eqnarray}
v_{1} &=&\mathbb{E}_{\left( X,Y\right) }\left[ \left( 2\left( Y-s_{M}\left(
X\right) \right) -s_{n_{1}}\left( M\right) \left( X\right) +s_{M}\left(
X\right) \right) ^{2}\left( s_{n_{1}}\left( M\right) \left( X\right)
-s_{M}\left( X\right) \right) ^{2}\right]  \notag \\
&\leq &\left( 4A+\left\Vert s_{n_{1}}\left( M\right) -s_{M}\right\Vert
_{\infty }\right) ^{2}\mathbb{E}_{X}\left[ \left( s_{n_{1}}\left( M\right)
\left( X\right) -s_{M}\left( X\right) \right) ^{2}\right]  \notag \\
&=&\left( 4A+\left\Vert s_{n_{1}}\left( M\right) -s_{M}\right\Vert _{\infty
}\right) ^{2}P\left( Ks_{n_{1}}\left( M\right) -Ks_{M}\right)  \label{v1}
\end{eqnarray}%
and 
\begin{eqnarray}
b_{1} &=&\left\Vert \left( 2\left( Y-s_{M}\left( X\right) \right)
-s_{n_{1}}\left( M\right) \left( X\right) +s_{M}\left( X\right) \right)
\left( s_{n_{1}}\left( M\right) \left( X\right) -s_{M}\left( X\right)
\right) \right\Vert _{\infty }  \notag \\
&\leq &4A\left\Vert s_{n_{1}}\left( M\right) -s_{M}\right\Vert _{\infty
}+\left\Vert s_{n_{1}}\left( M\right) -s_{M}\right\Vert _{\infty }^{2}\text{
.}  \label{b1}
\end{eqnarray}%
Now, we set $\Omega _{v}=\left\{ v_{1}\leq L_{v}\left( D_{M}\vee \ln
n_{1}\right) /n_{1}\right\} $ and $\Omega _{b}=\left\{ b_{1}\leq L_{b}\sqrt{%
D_{M}\ln n_{1}/n_{1}}\right\} $. By integrating (\ref{condi_ber}), it comes
for all $x>0$,%
\begin{eqnarray*}
&&\mathbb{P}\left( \left\vert P_{n_{2}}\left( Ks_{n_{1}}\left( M\right)
-Ks_{M}\right) -P\left( Ks_{n_{1}}\left( M\right) -Ks_{M}\right) \right\vert
\geq x\right) \\
&\leq &2\mathbb{E}\left[ \exp \left( -\frac{n_{2}x^{2}}{2\left(
v_{1}+b_{1}x/3\right) }\right) \mathbf{1}_{\Omega _{v}\cap \Omega _{b}}%
\right] +2\mathbb{P}\left( \Omega _{v}^{c}\right) +2\mathbb{P}\left( \Omega
_{b}^{c}\right) \\
&\leq &2\exp \left( -\frac{n_{2}x^{2}}{2\left( L_{v}\left( D_{M}\vee \ln
n_{1}\right) /n_{1}+L_{b}x\sqrt{D_{M}\ln n_{1}/n_{1}}\right) }\right) +2%
\mathbb{P}\left( \Omega _{v}^{c}\right) +2\mathbb{P}\left( \Omega
_{b}^{c}\right)
\end{eqnarray*}%
From assumption (\textbf{Ac}$_{\infty }$) and inequality (\ref%
{upper_true_selection}), it is possible to choose $L_{v}$ and $L_{b}$,
depending among other constants on $c$, such that for all $n\geq n_{0}\left( 
\text{(\textbf{GSA})},c\right) $, $2\mathbb{P}\left( \Omega _{v}^{c}\right)
+2\mathbb{P}\left( \Omega _{b}^{c}\right) \leq 10n^{-2-\alpha _{\mathcal{M}%
}} $. Thus, we get for $L>0$ large enough and for all $x>0$,%
\begin{eqnarray}
&&\mathbb{P}\left( \left\vert P_{n_{2}}\left( Ks_{n_{1}}\left( M\right)
-Ks_{M}\right) -P\left( Ks_{n_{1}}\left( M\right) -Ks_{M}\right) \right\vert
\geq x\right)  \notag \\
&\leq &2\exp \left( -\frac{n_{2}x^{2}}{L\left( \left( D_{M}\vee \ln
n_{1}\right) /n_{1}+x\sqrt{D_{M}\ln n_{1}/n_{1}}\right) }\right)
+10n^{-2-\alpha _{\mathcal{M}}}\text{ .}  \label{p1n2_control}
\end{eqnarray}%
By taking $x=\sqrt{L\alpha \ln n\left( D_{M}\vee \ln n_{1}\right) \left(
L\alpha \left( \ln n\right) \left( \ln n_{1}\right) +4n_{2}\right) }/\left(
n_{2}\sqrt{n_{1}}\right) >0$ in the latter inequality, it comes%
\begin{gather*}
\mathbb{P}\left( \left\vert P_{n_{2}}\left( Ks_{n_{1}}\left( M\right)
-Ks_{M}\right) -P\left( Ks_{n_{1}}\left( M\right) -Ks_{M}\right) \right\vert
\geq L\frac{\sqrt{\left( D_{M}\vee \ln n_{1}\right) \left( \ln n\right)
\left( \left( \ln n\right) \left( \ln n_{1}\right) +n_{2}\right) }}{n_{2}%
\sqrt{n_{1}}}\right) \\
\leq 12n^{-2-\alpha _{\mathcal{M}}}\text{ ,}
\end{gather*}%
where $L>0$ depends on the constants in (\textbf{GSA}) and on $c$.
Inequalities (\ref{p1n2_3}) and (\ref{p1n2_4}) then follow from simple
calculations.

\bigskip
\end{proof}

\begin{remark}
\label{remark_p1n2}It is easy to see that by using the assumption of
consistency in sup-norm for a fixed model, stated as (\textbf{H5}) in \cite%
{saum:12}, instead of (\textbf{Ac}$_{\infty }$) and by using Theorem 4\ of\ 
\cite{saum:12} instead of inequality (\ref{upper_true_selection}), the
results established in Lemma \ref{Lemma_p1n2} are valid with probability
bounds proportional to $n^{-\alpha }$, for any $\alpha >0$ (in Lemma \ref%
{Lemma_p1n2}, we only derive the case $\alpha =2+\alpha _{\mathcal{M}}$ for
convenience).
\end{remark}

\bigskip

\noindent \textbf{Proof of Theorem \ref{theorem_pen_n1n2}. }We set $%
\pen%
_{0}\left( M\right) =%
\pen%
_{ho}\left( M\right) -\left( n_{1}/n\right) \cdot \left( P_{n_{2}}\left(
Ks_{\ast }\right) -P_{n_{1}}\left( Ks_{\ast }\right) \right) $. It is worth
noting that $P_{n_{2}}\left( Ks_{\ast }\right) -P_{n_{1}}\left( Ks_{\ast
}\right) $ is a quantity independent of $M$, when $M$ varies in $\mathcal{M}%
_{n}$. Hence, the procedure defined by $%
\pen%
_{0}$ gives the same result as the hold-out procedure defined by $%
\pen%
_{ho}$. It will be convenient for our analysis to consider $%
\pen%
_{0}$ instead of $%
\pen%
_{ho}$. As a matter of fact, we derive Theorem \ref{theorem_pen_n1n2}\ as a
corollary of Theorem \ref{theorem_opt_pen_reg_pp}\ applied with $%
\pen%
\equiv 
\pen%
_{0}$, through the use of Lemma \ref{Lemma_p1n2}.

We get for all $M\in \mathcal{M}_{n}$,%
\begin{eqnarray*}
\pen%
_{0}\left( M\right) &=&\frac{n_{1}}{n}\left( P_{n_{2}}\left(
Ks_{n_{1}}\left( M\right) -Ks_{\ast }\right) -P_{n_{1}}\left(
Ks_{n_{1}}\left( M\right) -Ks_{\ast }\right) \right) \\
&=&\frac{n_{1}}{n}\left( P_{n_{2}}\left( Ks_{n_{1}}\left( M\right)
-Ks_{M}\right) -P_{n_{1}}\left( Ks_{n_{1}}\left( M\right) -Ks_{M}\right)
\right) \\
&&+\frac{n_{1}}{n}\left( \left( P_{n_{2}}-P\right) \left( Ks_{M}-Ks_{\ast
}\right) -\left( P_{n_{1}}-P\right) \left( Ks_{M}-Ks_{\ast }\right) \right)
\\
&=&\frac{n_{1}}{n}\left( \text{p}_{1}^{n_{2}}\left( M\right) +\text{p}%
_{2}^{n_{1}}\left( M\right) +\bar{\delta}^{n_{2}}\left( M\right) -\bar{\delta%
}^{n_{1}}\left( M\right) \right)
\end{eqnarray*}%
where%
\begin{equation*}
\text{p}_{1}^{n_{2}}\left( M\right) =P_{n_{2}}\left( Ks_{n_{1}}\left(
M\right) -Ks_{M}\right) \text{ , p}_{2}^{n_{1}}\left( M\right)
=P_{n_{1}}\left( Ks_{M}-Ks_{n_{1}}\left( M\right) \right) \text{ , }\bar{%
\delta}^{n_{i}}\left( M\right) =\left( P_{n_{i}}-P\right) \left(
Ks_{M}-Ks_{\ast }\right) \text{ .}
\end{equation*}%
Let $\Omega _{n}$ be the event on which:

\begin{itemize}
\item For all models $M\in \mathcal{M}_{n}$ of dimension $D_{M}$ such that $%
A_{\mathcal{M},+}\left( \ln n\right) ^{3}\leq D_{M}$, it holds 
\begin{align}
\left\vert 
\p_{1}\left(M\right)%
-\mathbb{E}\left[ \text{p}_{2}\left( M\right) \right] \right\vert & \leq L_{%
\text{(\textbf{GSA})}}\varepsilon _{n}\left( M\right) \mathbb{E}\left[ \text{%
p}_{2}\left( M\right) \right]  \label{line_1_n1} \\
\left\vert \text{p}_{2}\left( M\right) -\mathbb{E}\left[ \text{p}_{2}\left(
M\right) \right] \right\vert & \leq L_{\text{(\textbf{GSA})}}\varepsilon
_{n}^{2}\left( M\right) \mathbb{E}\left[ \text{p}_{2}\left( M\right) \right]
\label{line_2_n1}
\end{align}%
together with%
\begin{eqnarray}
\left\vert \text{p}_{1}^{n_{2}}\left( M\right) -\frac{n}{n_{1}}\mathbb{E}%
\left[ \text{p}_{2}\left( M\right) \right] \right\vert &\leq &L_{\text{(%
\textbf{GSA})},c}\left[ \varepsilon _{n}^{1,2}\left( M\right) +\varepsilon
_{n}\left( M\right) \right] \mathbb{E}\left[ \text{p}_{2}\left( M\right) %
\right]  \label{line_p1n2} \\
\left\vert \text{p}_{2}^{n_{1}}\left( M\right) -\frac{n}{n_{1}}\mathbb{E}%
\left[ \text{p}_{2}\left( M\right) \right] \right\vert &\leq &L_{\text{(%
\textbf{GSA})},c}\varepsilon _{n}^{2}\left( M\right) \mathbb{E}\left[ \text{p%
}_{2}\left( M\right) \right]  \label{line_p2n1} \\
\left\vert \bar{\delta}^{n_{1}}\left( M\right) \right\vert &\leq &\frac{\ell
\left( s_{\ast },s_{M}\right) }{\sqrt{D_{M}}}+L_{\text{(\textbf{GSA})},c}%
\frac{\ln n}{\sqrt{D_{M}}}\mathbb{E}\left[ \text{p}_{2}\left( M\right) %
\right]  \label{line_delta_n1_bis} \\
\left\vert \bar{\delta}^{n_{2}}\left( M\right) \right\vert &\leq &L_{\text{(%
\textbf{GSA})}}\left( \sqrt{\frac{\ell \left( s_{\ast },s_{M}\right) \ln
n_{2}}{n_{2}}}+\frac{\ln n_{2}}{n_{2}}\right)  \label{line_delta_n2_bis}
\end{eqnarray}

\item For all models $M\in \mathcal{M}_{n}$ of dimension $D_{M}$ such that $%
D_{M}\leq A_{\mathcal{M},+}\left( \ln n\right) ^{3}$, it holds%
\begin{eqnarray}
\left\vert \bar{\delta}^{n_{1}}\left( M\right) \right\vert &\leq &L_{\text{(%
\textbf{GSA})},c}\left( \sqrt{\frac{\ell \left( s_{\ast },s_{M}\right) \ln n%
}{n}}+\frac{\ln n}{n}\right)  \label{delta_n1_small} \\
\left\vert \bar{\delta}^{n_{2}}\left( M\right) \right\vert &\leq &L_{\text{(%
\textbf{GSA})},c}\left( \sqrt{\frac{\ell \left( s_{\ast },s_{M}\right) \ln n%
}{n}}+\frac{\ln n}{n}\right)  \label{delta_n2_small} \\
\text{p}_{2}^{n_{1}}\left( M\right) &\leq &L_{\text{(\textbf{GSA})},c}\frac{%
D_{M}\vee \ln n}{n}\leq L_{\text{(\textbf{GSA})},c}\frac{\left( \ln n\right)
^{3}}{n}  \label{p2n1_small} \\
\text{p}_{1}^{n_{2}}\left( M\right) &\leq &L_{\text{(\textbf{GSA})},c}\left( 
\frac{\left( \ln n\right) ^{2}}{n}+\frac{D_{M}\vee \ln n}{n}\right) \leq L_{%
\text{(\textbf{GSA})},c}\frac{\left( \ln n\right) ^{3}}{n}
\label{p1n2_small}
\end{eqnarray}
\end{itemize}

\noindent By (\ref{p1_selection}), (\ref{p2_selection}), (\ref%
{upper_true_selection}) and (\ref{upper_emp_selection}) in Remark \ref%
{remark_application_fixed_model}, Lemma \ref{mean_emp_risk} and Lemma \ref%
{Lemma_p1n2}, we get for all $n\geq n_{0}\left( \text{(\textbf{GSA})}%
,c\right) $, 
\begin{equation*}
\mathbb{P}\left( \Omega _{n}\right) \geq 1-A_{p}n^{-2}-L\sum_{M\in \mathcal{M%
}_{n}}n^{-2-\alpha _{\mathcal{M}}}\geq 1-L_{A_{p},c_{\mathcal{M}}}n^{-2}%
\text{ }.
\end{equation*}

\bigskip

\noindent We consider models $M\in \mathcal{M}_{n}$ such that $A_{\mathcal{M}%
,+}\left( \ln n\right) ^{3}\leq D_{M}$. Notice that (\ref{line_delta_n1_bis}%
) implies by (\ref{def_epsilon_selection}) that, for all $M\in \mathcal{M}%
_{n}$ such that $A_{\mathcal{M},+}\left( \ln n\right) ^{3}\leq D_{M}$, 
\begin{eqnarray*}
\left\vert \bar{\delta}^{n_{1}}\left( M\right) \right\vert &\leq &L_{\text{(%
\textbf{GSA})},c}\left( \frac{\left( \ln n\right) ^{3}}{D_{M}}%
\cdot%
\frac{\ln n}{D_{M}}\right) ^{1/4}\times \left( \ell \left( s_{\ast
},s_{M}\right) +\mathbb{E}\left[ p_{2}\left( M\right) \right] \right) \\
&\leq &L_{\text{(\textbf{GSA})},c}\varepsilon _{n}\left( M\right) \left(
\ell \left( s_{\ast },s_{M}\right) +\mathbb{E}\left[ p_{2}\left( M\right) %
\right] \right) \text{ .}
\end{eqnarray*}%
In addition, from (\ref{line_delta_n2_bis}), Lemma \ref{mean_emp_risk} and
the fact that $n\left( \ln n\right) ^{\tau }/D_{M}\leq n_{2}$, we get that
for all $n\geq n_{0}\left( \text{\textbf{GSA}}\right) $,%
\begin{eqnarray*}
\left\vert \bar{\delta}^{n_{2}}\left( M\right) \right\vert &\leq &L_{\text{(%
\textbf{GSA})}}\left( \sqrt{\frac{\ell \left( s_{\ast },s_{M}\right) \ln
n_{2}}{n_{2}}}+\frac{\ln n_{2}}{n_{2}}\right) \\
&\leq &L_{\text{(\textbf{GSA})}}\left( \frac{\ell \left( s_{\ast
},s_{M}\right) }{\left( \ln n\right) ^{\left( \tau -1\right) /2}}+\frac{\ln
n_{2}}{n_{2}}\left( \ln n\right) ^{\left( \tau -1\right) /2}\right) \\
&\leq &L_{\text{(\textbf{GSA})}}\left( \ln n\right) ^{\left( 1-\tau \right)
/2}\left( \ell \left( s_{\ast },s_{M}\right) +\mathbb{E}\left[ p_{2}\left(
M\right) \right] \right) \text{ .}
\end{eqnarray*}%
We deduce that on $\Omega _{n}$ we have, for all models $M\in \mathcal{M}%
_{n} $ such that $A_{\mathcal{M},+}\left( \ln n\right) ^{3}\leq D_{M}$ and
for all $n\geq n_{0}\left( \text{\textbf{GSA}}\right) $,

\begin{eqnarray}
&&\left\vert 
\pen%
_{0}\left( M\right) -2\mathbb{E}\left[ \text{p}_{2}\left( M\right) \right]
\right\vert  \notag \\
&\leq &\frac{n_{1}}{n}\left( \left\vert \text{p}_{1}^{n_{2}}\left( M\right) -%
\frac{n}{n_{1}}\mathbb{E}\left[ \text{p}_{2}\left( M\right) \right]
\right\vert +\left\vert \text{p}_{2}^{n_{1}}\left( M\right) -\frac{n}{n_{1}}%
\mathbb{E}\left[ \text{p}_{2}\left( M\right) \right] \right\vert \right) 
\notag \\
&&+\left\vert \bar{\delta}^{n_{1}}\left( M\right) \right\vert +\left\vert 
\bar{\delta}^{n_{2}}\left( M\right) \right\vert  \notag \\
&\leq &\left( L_{\text{(\textbf{GSA})},c}\left( \varepsilon _{n}^{1,2}\left(
M\right) +\varepsilon _{n}\left( M\right) +\left( \ln n\right) ^{\left(
1-\tau \right) /2}\right) \right) \left( \ell \left( s_{\ast },s_{M}\right) +%
\mathbb{E}\left[ p_{2}\left( M\right) \right] \right)
\label{control_pen_ho_res}
\end{eqnarray}%
Hence, inequality (\ref{pen_id_pp}) of Theorem \ref{theorem_opt_pen_reg_pp}\
is satisfied on $\Omega _{n}$ by taking 
\begin{equation*}
\delta =L_{\text{(\textbf{GSA})},c}\left( \varepsilon _{n}^{1,2}\left(
M\right) +\varepsilon _{n}\left( M\right) +\left( \ln n\right) ^{\left(
1-\tau \right) /2}\right) \text{ .}
\end{equation*}%
Moreover, we have $\delta \in \left[ 0,1\right) $ for all $n\geq n_{0}\left( 
\text{(\textbf{GSA}),}c,\tau \right) $.

\noindent Let us now consider models $M\in \mathcal{M}_{n}$ such that $%
D_{M}\leq A_{\mathcal{M},+}\left( \ln n\right) ^{3}$. By (\ref%
{delta_n1_small}), (\ref{delta_n2_small}), (\ref{p1n2_small}) and (\ref%
{p2n1_small}), we have on $\Omega _{n}$,%
\begin{eqnarray}
\left\vert 
\pen%
_{0}\left( M\right) \right\vert &=&\frac{n_{1}}{n}\left\vert \text{p}%
_{1}^{n_{2}}\left( M\right) +\text{p}_{2}^{n_{1}}\left( M\right) +\bar{\delta%
}^{n_{2}}\left( M\right) -\bar{\delta}^{n_{1}}\left( M\right) \right\vert 
\notag \\
&\leq &L_{\text{(\textbf{GSA})},c}\left( \sqrt{\frac{\ell \left( s_{\ast
},s_{M}\right) \ln n}{n}}+\frac{\left( \ln n\right) ^{3}}{n}\right)  \notag
\\
&\leq &L_{\text{(\textbf{GSA})},c}\left( \frac{\ell \left( s_{\ast
},s_{M}\right) }{\left( \ln n\right) ^{2}}+\frac{\left( \ln n\right) ^{3}}{n}%
\right)  \label{control_pen_ho_small}
\end{eqnarray}%
Inequality (\ref{control_pen_ho_small}) implies that inequality (\ref%
{pen_id_2_pp}) of Theorem \ref{theorem_opt_pen_reg_pp}\ is satisfied with $%
A_{r}=L_{\text{(\textbf{GSA})},c}$. From (\ref{control_pen_ho_res}) and (\ref%
{control_pen_ho_small}), we thus apply Theorem \ref{theorem_opt_pen_reg_pp}\
with $A_{p}=L_{A_{p},c_{\mathcal{M}}}$, and this gives Theorem \ref%
{theorem_pen_n1n2} with%
\begin{equation*}
\theta _{n}=L_{\text{(\textbf{GSA})},c}\left( \left( \ln n\right)
^{-2}+\left( \ln n\right) ^{\left( 1-\tau \right) /2}+\sup_{M\in \mathcal{M}%
_{n}}\left\{ \varepsilon _{n}\left( M\right) +\varepsilon _{n}^{1,2}\left(
M\right) ,\text{ }A_{\mathcal{M},+}\left( \ln n\right) ^{3}\leq D_{M}\leq
n^{\eta +1/\left( 1+\beta _{+}\right) }\right\} \right) \text{ .}
\end{equation*}

\noindent {\LARGE Acknowledgements}\bigskip

I am deeply grateful to Pr. Jon A. Wellner and Pr. Pascal Massart for their
valuable support. I also warmly thank Pr. Wellner for having helped me to
improve my English along the text. Finally, I gratefully thank the associate
editors and anonymous referees for their comments and suggestions, that
greatly improved the quality of the paper.

\bibliographystyle{plain}
\bibliography{Slope_heuristics_regression_13}

\end{document}